\magnification=1200
\hsize=16 true cm    \vsize=22.7 true cm
\hoffset=4 true mm   \voffset=2 true mm

\font\tengoth = eufm10
\font\sevengoth = eufm7
\font\fivegoth = eufm5

\newfam\gothfam
\textfont\gothfam=\tengoth
\scriptfont\gothfam=\sevengoth
\scriptscriptfont\gothfam=\fivegoth

\def\goth{\tengoth\fam\gothfam}

\font\tenmath =msbm10
\font\sevenmath=msbm7

\newfam\mathfam
\textfont\mathfam=\tenmath
\scriptfont\mathfam=\sevenmath
\scriptscriptfont\mathfam=\sevenmath
\def\math{\tenmath\fam\mathfam}
\def\build#1_#2{\mathrel{\mathop{\kern 0pt#1}\limits_{#2}}}
\def\K{{\math K}}
\def\wt{\hbox{wt}}\def\Ker{\hbox
{Ker}}
\def\dim{\hbox{dim}}
\def\Gr{\hbox{Gr}}%
\def\ht{\hbox{ht}}
\def\hb{{\math H}_q}
\def\dim{\hbox{dim}}\def\ad{\hbox{ad}}
\def\K{{\math K}}
\def\sl{\hbox{{\goth sl}}}
\def\build#1_#2{\mathrel{\mathop{\kern 0pt#1}\limits_{#2}}}
\def\B{{\cal B}}\def\A{{A}}\def\E{{\cal E}}

\def\qed{\hfill$\diamondsuit$}
\def\M{{\cal M}}\def\co{{\cal C}_0}

\def\ua{U_{\A}}\def\uao{\ua^0}
\def\uap{\ua^+}\def\uam{\ua^-}\def\um{U(\nm)}

\def\R{{\cal R}}
\def\uq{U_q(\g)}\def\up{U_q(\n)}\def\um{U_q(\nm)}\def\uo{U_q^0}\def\uqo{U_q^0}
\def\lq{V_q}\def\lql{\lq(\lambda)}

\def\za{Z_{\A}}\def\ja{J_{\A}}\def\ha{H_{\A}}
\def\hf{\hat F_q}


\def\hb{{\math H}_q}\def\hq{\hb}


\def\N{{\math N}}\def\Z{{\math Z}}%
\def\Q{{\math Q}}%
\def\n{{\goth n}}%
\def\m{\hbox{{\goth m}}^+}%
\def\C{{\math C}}%
\def\g{{\goth g}}%
\def\h{{\goth h}}%
\def\b{{\goth b}}\def\nm{{\goth n}^-}%
\def\fa{F(\ua)}\def\far{F_{\A,\lambda}}\def\faru{F_{\A,\mu}}
\def\fua{F(\ua)}


\def\P{{\cal P}}
\def\Q{{\bf Q}}\def\rp{\Delta^+}

\font\petittitre=cmbx8 scaled 1000 
\def\H{{\cal H}}

\vskip 2cm\centerline{\petittitre{ ON HARMONIC ELEMENTS FOR SEMI-SIMPLE
LIE ALGEBRAS}}
\centerline{\petittitre{ }}
\vskip 1cm\centerline{Philippe CALDERO}
\vskip 4cm\noindent
\par\noindent{\petittitre ABSTRACT.} {\it Let $\g$ be a semi-simple
complex Lie
algebra, $\g=\nm\oplus\h\oplus\n$ its triangular decomposition. Let
$U(\g)$, resp. $\uq$, be its enveloping algebra, resp. its quantized
enveloping
algebra.
 This article gives a quantum approach
to the combinatorics of (classical)
harmonic elements and Kostant's generalized exponents for $\g$.
A quantum analogue of
the space of harmonic elements has been given in [21]. On the one hand,
we give specialization results concerning harmonic elements, central
elements of $\uq$, and the
Joseph and Letzter's decomposition, see [21]. For
$\g=\sl_{n+1}$, we describe the specialization
of quantum harmonic space in the $\N$-filtered algebra
$U(\sl_{n+1})$ as the materialization of a theorem of
Lascoux-Leclerc-Thibon, [29].
This enables us to study
a Joseph-Letzter decomposition in the algebra
$U(\sl_{n+1})$. On the other hand, we prove that highest weight 
harmonic elements
 can be  calculated
in terms of the dual of Lusztig's canonical base.
In the simply laced case, we parametrize  a base of $\n$-invariants of
minimal primitive
quotients by the set $\co$ of integral points of a convex cone.}
  \vskip 2cm\noindent
This research has been partially supported by the EC TMR network
 "Algebraic Lie Representations" , contract no. ERB FMTX-CT97-0100.

\vskip 1cm\noindent
\noindent
{\bf 0. Introduction.}\bigskip\noindent
{\bf 0.1.} Let $V$ be a complex finite dimensional space and $G$ a
Lie subgroup of
GL$(V)$. Assume that $V$ is completely reducible as a $G$-module.
 Then $G$ acts semi-simply on the symmetric $\C$-algebra $S(V)$.  Let
$J$ be the ideal generated by the non-constant homogeneous $G$-invariant
elements of $S(V)$. Let $V^*$ be the dual of $V$.
There exists a natural non-degenerate duality between $S(V)$
and $S(V^*)$. Let $H^*\subset S(V^*)$ be the orthogonal of $J$, $H^*$ is
the space of $G$-harmonic polynomials of $V$. \par
Suppose now that
$G$ is a complex simply connected semi-simple Lie group acting by the
adjoint
action on its associated Lie algebra $\g=V$. Then, the Killing form
maps the harmonic space $H^*$  onto $H\subset S(\g)$.
We have $S(\g)=H\oplus J$, cf [27]. \par
Let $U(\g)$ be the enveloping algebra of $\g$ and $Z(\g)$ be its center.
Recall that $U(\g)$  is naturally filtred and the associated graded
algebra of 
$U(\g)$  is $S(\g)$,  by the Poincar\'e-Birkhoff-Witt
theorem. As in [27], we can choose a sub-ad$\g$-module
 $\H$ of $U(\g)$ such the graded space associated to $\H$ is $H$.
The space $\H$ will be called
harmonic space in $U(\g)$. Kostant's separation theorem states that
$U(\g)\simeq Z(\g)\otimes {\H}$ via multiplication.\smallskip\noindent
{\bf 0.2.} Fix a Cartan subalgebra $\h$ of $\g$ and let
$\g=\nm\oplus\h\oplus\n$ be the triangular decomposition of $\g$.
Let $P$ be the group of integral weights of $\g$ and $P^+$ be the
semi-group of dominant weights of $\g$, generated by the fundamental
weights $\varpi_i$, $1\leq i\leq n$. The set
$\Delta$, resp. $\Delta^+$, is the root system, resp. the set of positive
roots,
and $W$ is the Weyl group. For all $\lambda$ in $P^+$, let $V(\lambda)$ be
the simple $\g$-module with highest weight $\lambda$. \par
Let's sketch a few results on ${H}$ that can be found in [27].
First of all, ${H}$ is a $\ad\g$-module and the multiplicity of
$V(\lambda)$ in $H$ is
$$k_\lambda:=[{H}\, :\,V(\lambda)]=\dim\,V(\lambda)_0,\leqno (0.2.1)$$
 where  $V(\lambda)_0$ is the
0-weight subspace of $V(\lambda)$. Let $H(\lambda)$ be the isotypical
component
of $V(\lambda)$ in $H$ and let $S(\g)^n$ be the
$n$-th graded component of $S(\g)$.
Then, there are integers $m_\lambda^1\leq\ldots m_\lambda^{k_\lambda}$
such that
$$\sum_{n\geq 0}[H(\lambda)\cap S(\g)^n \,:\,V(\lambda)]t^n =
\sum_{i=1}^{k_\lambda}t^{m_\lambda^i},\leqno (0.2.2)$$
as polynomials in $t$.
The  $m_\lambda^i$ are called Kostant's generalized exponents. Indeed,
they
generalize the usual exponents of $\g$ because they can be realized as
exponents of the eigenvalue of the Coxeter element acting on the
zero weight space of $V(\lambda)$. Two problems
arise :\smallskip\noindent
{\bf Problem 1.} Find some combinatorial rules that calculates the
generalized
exponents.
\smallskip\noindent
{\bf Problem 2.} Calculate the harmonic elements of
$U(\g)$.\smallskip\noindent
{\bf 0.3.} Let's sketch results about Problem 1. For $\mu\in P$, let
$\P(\mu)$ be the dimension of the space of weight $\mu$ in $S(\n)$. $\P$
is
 Kostant's partition function. In the completed group algebra $\C<P>$,
cf. [14, 7.5],
 we have $\sum_{\mu\in P}\P(\mu)e^{-\mu}=\prod_{\alpha\in \Delta^+}
(1-e^{-\alpha})^{-1}.$ Lusztig's $q$-analogue of Kostant's partition
function
is given by $\sum_{\mu\in P}\P_q(\mu)e^{-\mu}=\prod_{\alpha\in \Delta^+}
(1-qe^{-\alpha})^{-1},$ where $q$ is an indeterminate.
For $\lambda\in P^+$ and $\mu\in P$,
the Kostka-Foulkes polynomial is defined by $$K_{\lambda\mu}(q)
=\sum_{w\in W}\,(-1)^{l(w)}\P_q(w(\lambda+\rho)-\mu-\rho),\leqno (0.3.1)$$
where
$l(w)$ is the length of $w$ and $\rho$ is the sum of the fundamental
weights.
By [14, 7.5.10], $K_ {\lambda\mu}(1)=\dim\,V(\lambda)_\mu$.
R. Brylinski gives the following interpretation of Hesselink formula :
\medskip\noindent
{\bf Theorem.}[17], [7] {\it For $\lambda\in P^+$,
$K_{\lambda\,0}(q)=\sum_{i=1}^{k_\lambda}\,
q^{m_\lambda^i}$.}\medskip\noindent
But this is not an answer to Problem 1 because of the signs in (0.3.1).
\smallskip
In [29], there is given some nice combinatorics for the generalized
exponents in the case A$_n$ : first a multi-variable Kostka-Foulkes
function ${\K}_ {\lambda\,0}(X_1,\ldots,X_n)$ is defined
in terms of the set $B(\lambda)_0$
of zero weight elements in the crystal base of $V(\lambda)$ by
$${\K}_ {\lambda\,0}(X_1,\ldots,X_n)=\sum_{b\in B(\lambda)_0}
\prod_{i=1}^n\,X_i^{\varepsilon_i(b)},$$
where $\varepsilon_i$ are the usual parameters of the crystal base.
It can be easily
calculated by the combinatorics described in [25]. Then, the
Kostka-Foulkes
function $K_{\lambda\,0}(q)$ is recovered by specializing
$X_i$ on $q^i$ :\medskip\noindent
{\bf Theorem.} [29] {\it For $\lambda\in P^+$,
$K_{\lambda\,0}(q)={\K}_
{\lambda\,0}(q,q^2,,\ldots,q^n)$.}\medskip\noindent
So much for Problem 1 in the A$_n$ case. Now, Problem 2 has no known
general solution even in $\sl_3$.
\medskip\noindent
{\bf 0.4.} The article aims to present a quantum theory
of harmonic spaces for enveloping algebras of semi-simple Lie algebras,
in order to use it to get results in the classical case.\par
Let $q$ be an indeterminate. Let $\uq$ be the
 $q$-deformation of the enveloping algebra $U(\g)$. For
$\lambda$ in $P^+$, let $V_q(\lambda)$ be the quantum simple module with
highest
weight $\lambda$.\par
The algebra $\uq$ contains
the algebra $F(\uq)$ of ad-finite elements of  $\uq$, see 1.4. We use
Joseph-Letzter's decomposition theorem
as a basic tool, [21]. It asserts that  $$F(\uq)=\bigoplus_{\lambda\in
P^+}
F_q^\lambda,\leqno (0.4.1)$$ where $F_q^\lambda$ are sub-ad$\uq$-modules
of 
$F(\uq)$ isomorphic to
$ V_q(\lambda)^*\otimes
V_q(\lambda)$
endowed with the diagonal action. This can be seen, [9],
as a Peter-Weyl theorem embedded in $\uq$ via the Rosso form, [34].
The center $Z$ of $\uq$ is as in the classical case a polynomial algebra
with
dimension $n=rk(\g)$. The invariant elements $z_i$ of $F_q^{\varpi_i}$,
where
$1\leq i\leq n$, are algebraically independent generators of $Z$.
Let $J_q$ be the ideal in $F(\uq)$
generated by the
$z_i$. By definition, an harmonic space $H_q$ will be a complementary 
$\ad\uq$-module of $J_q$
in $F(\uq)$, such that $H_q=\oplus_{\lambda\in P^+}F_q^\lambda\cap H_q$,
with compatibility conditions, see [1.8, Definition 1].
\medskip\noindent
{\bf 0.5.} The second tool is what we call Ringel's filtration :
a Poincar\'e-Birkhoff-Witt base of $\uq$ can be lexicographically
ordered in such a way that
the  graded algebra associated to this ordering, see 1.3, is
$q$-commutative.
We prove that, for this ordering :\smallskip\noindent
{\bf Theorem 1.} {\it The set
$\{z_i,\,1\leq i\leq n\}$ is a Gr\"obner base of the ideal
$J_q$.}\smallskip\noindent
 We obtain that $\uq$
is free over its center, with a set of explicit generators given in terms
of
roots packages, [11]. We prove an analogue of this result for quantized
enveloping algebras at a root of one.\smallskip\noindent
{\bf Theorem 2.} {\it When $q=\varepsilon$ is a root of one, then the
algebra $U_{\varepsilon}$, see [12, 1.5], is free over its
center.}\smallskip\noindent
As before, an explicit set of generators can be given.
\medskip\noindent
{\bf 0.6.} For $\lambda\in P^+$, let $H_q(\lambda)$ be the isotypical
component
of type $V_q(\lambda)$ in $H_q$. A result of Joseph and Letzter asserts
that
(0.2.1) holds in the quantum case. Now, Joseph-Letzter's decomposition,
see (0.4.1), affords a natural $P^+$-filtration on $F(\uq)$.
By analogy with equation (0.2.2), we define
$P^+$-exponents of $\lambda$ as the family $\mu_\lambda^1$,$\ldots$,
$\mu_\lambda^{k_\lambda}$ of dominant weights such that
$$\sum_{\mu\in P^+}[H_q(\lambda)\cap \ad\uq.K_{-2\mu}
\,:\,V_q(\lambda)]e^\mu =
\sum_{i=1}^{k_\lambda}e^{\mu_\lambda^i}$$
in the group algebra of $P$.
 We prove, see Proposition 1.8
and Proposition 3.4, that the quantum analogues of
Problem 1 and Problem 2 can be solved in terms of
Lusztig's canonical base and its dual : \smallskip\noindent
{\bf Proposition.} {\it Let $\mu$ in $P^+$.
Let $\B_H$ be the set of elements $b$
of the  canonical base such that $\varepsilon_i(b)\leq -(\wt(b)
,\alpha_i\check{\,})$, $1\leq i\leq n$. Then,\par\noindent
\item{(i)} $\sum_\lambda[H_q\cap F_q^\lambda\,:V_q(\mu)]\prod_i
X_i^{<\lambda,\alpha_i
 \check{\,}>}
={\K}_ {\mu\,0}(X_1,\ldots,X_n)$,\par\noindent
\item{(ii)} a base of the space of $\n$-invariants of $H_q$ is labelled
by $\B_H$.}\smallskip\noindent
Note that (i) can be also easily deduced from a result of Baumann, [1].
In particular, for $\g=\sl_{n+1}$,
the multi-variable
Kostka-Foulkes function ${\K}_ {\lambda\,0}(X_1,\ldots,X_n)$
can be described as the Poincar\'e polynomial of $H_q(\lambda)$ in the
$P^+$-Joseph-Letzter
filtration.
\medskip\noindent
{\bf 0.7.} Another related problem is the following :\smallskip\noindent
{\bf Problem 3.} Describe a Joseph-Letzter decomposition inside the
classical
enveloping algebra $U(\g)$ endowed with its natural
filtration.\smallskip\noindent
This  includes the problem of specialization of harmonic elements as well
as the specialization of the center. For general $\g$, we give some 
preparation theorems. Let $A$ be
the algebra $\C[q]$ localized at $q=1$. We present a $A$-form $\ha$ which 
links the quantum harmonic space with the classical one. Indeed, it verifies
$H_q=\C(q)\otimes_A\ha$ and $H_1=\C\otimes_A\ha$ which both verify the separation theorem
at, respectively, the quantum and the classical level. We prove a
separation theorem on $A$-forms.\par
For $\g=\sl_{n+1},$ the understanding of the specialization of
a) the quantum harmonic space and b) the center, see [8], enables
us, by the separation theorem, to give an answer to Problem 3,
see Theorem 4.3. 
\medskip\noindent
{\bf 0.8.}
Recall that any reduced decomposition of the longest element $w_0\in W$ in
the
Weyl group gives rise to a Poincar\'e-Birkhoff-Witt basis of $\up$, whose
elements are
labelled by points in $\N^N$, where $N=$dim$\n$. When $\g$ has type
 A-D-E, then by [31], [32], the dual of the canonical
 base is in bijection with the PBW basis, with nice multiplicative
properties,
see (3.5.1)-(3.5.3). Hence, each element of the canonical base
 corresponds to an integral point in ${\math R}^N$, where
$N=\dim\n$.\smallskip\noindent
 {\bf Theorem.} {\it Let $\co\subset \N^N$ be the image of $\B_H$
 via the correspondence above. Then, $\co$ is the set of integral points
 in a convex cone.}\smallskip\noindent
 We describe a degenerescence of the variety $\C[\co]$ which relates the
 combinatorics of Point 2, see Theorem 3.5.
\bigskip\noindent
{\bf 1. Preliminaries and notations.}\bigskip\noindent
{\bf 1.1.} Let $\g$ be a semi-simple Lie $\C$-algebra  of rank $n$.
We fix a Cartan subalgebra   $\h$ of $\g$. Let $\g=\nm\oplus\h\oplus\n$
be the triangular decomposition  and $\{\alpha_i\}_i$ be a base of the
root
system
$\Delta$ resulting from this decomposition. We note $\b=\n+\h$ and
$\b^-=\nm+\h$ the two opposite Borel sub-algebras.
Let  $\rp$ be the set of positive roots, $P$ be the weight lattice
generated by the fundamental weights
$\varpi_i$, $1\leq i\leq n$,
and $P^+:=\sum_i\,\N\varpi_i$ the semigroup of integral dominant weights.
$P$ is endowed with the ordering $\preceq$
$$\lambda\preceq\mu\Leftrightarrow \mu-\lambda\in P^+.$$
We fix a {\it total} additive ordering $\leq$ on $P$ such that
$\mu-\lambda\in Q^+\Rightarrow \lambda\leq\mu \,\,\lambda,\mu\in P.$
Such an ordering always exists. For example, for $\lambda=\sum_i\lambda_i
\alpha_i$, $\mu=\sum_i\mu_i
\alpha_i$, $\lambda\leq\mu$ iff $(\sum_i\lambda_i,\lambda_1,\ldots,\lambda_n)$
preceeds $(\sum_i\mu_i,\mu_1,\ldots,\mu_n) $ for the lexicographical ordering
of ${\math Q}^{n+1}$ verifies the hypothesis.\par\noindent
 Let $W$ be the  Weyl group, generated by the
 reflections corresponding to the simple roots $s_i:=s_{\alpha_i}$. Let
$w_0$
be the longest element of $W$.
We note  $(\,,\,)$ a  $W$-invariant form on $P$. For $\beta\in\rp$,
$\ht(\beta)$
will mean the height of $\beta$.\medskip\noindent
{\bf 1.2.} Let $q$ be an indeterminate and $\uq$ be the simply connected
quantized
enveloping  algebra, defined as in [13, 0.2-0.3]. Let $\up$, resp. $\um$,
be the subalgebra generated by the canonical generators  $E_{\alpha_i}$,
resp. $F_{\alpha_i}$,
 of positive,  resp. negative, weights.
For all $\lambda$
in $P$, let $K_\lambda$ be the corresponding  element in the algebra
$\uo=\C(q)[P]$
of the torus of $\uq$. As in the classical case, we have the triangular
decomposition : $\uq\simeq\um\otimes\uo\otimes\up$.\smallskip\noindent
$\uq$ is endowed with a  structure of Hopf algebra with comultiplication
$\Delta$,
 antipode $S$ and augmentation $\varepsilon$, [18, 3.2.9].\par
We define in $\uq$ the left adjoint action by
$\ad\,v.u=v_{(1)}uS(v_{(2)}),$
where $\Delta(v)=v_{(1)}\otimes v_{(2)}$ with the Sweedler notations.\par
Let $M$ be a $\uq$-module and $\mu\in P$. Let $M_\mu$ be the space of
elements of weight $\mu$, i.e  $M_\mu:=
\{u\in M, K_{\alpha}m=q^{(\alpha,\mu)}m, \forall \alpha\in
P\}$.\par\noindent
\medskip\noindent
{\bf 1.3.}  We fix a decomposition of the longest element of the Weyl
group
$w_0=s_{i_1}\ldots s_{i_N}$, where $N=\dim\,\n$. Set
$\beta_k=s_{i_1}\ldots s_{i_{k-1}}(\alpha_{i_k})$, $1\leq k\leq N$.
We endow with an order the set $\rp$ of positive roots  :
$\beta_N<
\ldots\beta_2<\beta_1=\alpha_{i_1}$, see [12, 1.7].
For all $\beta$ in $\rp$, let  $E_{\beta}$, resp. $F_{\beta}$,
be the root elements of $\up$, resp. $\um$. Set $E_i=E_{\alpha_i}$,
$F_i=F_{\alpha_i}$, for each simple root.
For each $\psi$ in $\N^N$, set $E_{\psi}:=\prod_{l=N}^1
E_{\beta_l}^{\psi_l}$,
resp. $F_{\psi}:=\prod_{l=N}^1 F_{\beta_l}^{\psi_l}$. \par\noindent
 For
$\Gamma=(\phi,\lambda,\psi)\in
\N^N\times P\times\N^N$, set $X_{\Gamma}:=F_{\phi}K_{-\lambda}E_{\psi}$.
We know, cf. [12, 1.6],
that the $E_{\psi}$, resp. $F_{\phi}$, resp. $X_{\Gamma}$, form a
Poincar\'e-Birkhoff-Witt
base of $\up$, resp $\um$, resp. $\uq$, for $\psi\in\N^N$, resp.
$\phi\in\N^N$, resp. $\Gamma\in\N^N \times P\times\N^N$. To $\Gamma=
(\phi,\lambda,\psi)\in
\N^N\times P\times\N^N$, we associate $\tilde\Gamma=
(\sum (\phi_l+\psi_l)\ht(\beta_l),\phi,\psi)
\in
\N^{2N+1}$.
The set of subspaces
$$\{\bigoplus_{\tilde\Gamma'\leq\tilde\Gamma}\C(q)X_{\Gamma'}, \Gamma\in
\N^N\times P\times\N^N\},$$
for the lexicographic ordering on $\N^{2N+1}$ define a filtration of
algebra on $\uq$.
 In the sequel, this will be called the Ringel filtration,
[33].
The associated graded $\Gr\uq$ is generated by the
$\Gr F_\beta$, $\Gr K_\lambda$, $\Gr E_\beta$,
where $\beta\in\rp$, $\lambda\in P$, and $q$-commuting relations,
[12, Proposition 1.7]. In particular, we have
$$\Gr E_\alpha\Gr E_\beta=q^{(\alpha,\beta)}\Gr E_\beta\Gr
E_\alpha,\;\hbox{if}
\,\beta<\alpha\leqno (1.3.1)$$

\medskip\noindent
{\bf 1.4.}   For all $\lambda$ in $P^+$,
let $\lql$ be the simple $\uq$-module
with highest weight $\lambda$ and highest weight vector $v_\lambda$. Let
$\lql^*$ be its dual, endowed with
a structure of left
$\uq$-module twisted by the antipode. We define the ad-finite part of
$\uq$ by $$F(\uq):=\{u\in\uq, \dim\,\ad\,\uq(u)<+\infty\}.$$ We have the
theorem,
[18, 7.1], [9] :\medskip\noindent
{\bf Theorem.} {\it $F(\uq)$ has the following properties :\par\noindent
\item{(i)} $F(\uq)$ is an algebra and a semi-simple module for
the adjoint action,\par\noindent
\item{(ii)} $F(\uq)=\bigoplus_{\lambda\in P^+}\ad\uq.K_{-2\lambda}$,
\item{(iii)} The module $\ad\uq.K_{-2\lambda}$ is isomorphic to
$\lql^*\otimes\lql$.}\qed\medskip\noindent
{\bf 1.5.} As in [11, 2.1], we define the root packages
from a decomposition of $w_0$.
\medskip\noindent
{\bf Definition and notation.}
Fix $w_0=s_{i_1}\ldots s_{i_N}$. For $1\leq j\leq n$, we call root
packages
the sets
 $\rp_j:=\{\beta_l, i_l=j\}$. For $m$, $1\leq m\leq k:=$Card$\rp_j$,
we define $\alpha_{j,m}$ to be the $m$-th
element in the decreasing sequence of the roots of $\rp_j$ :
$\alpha_{j,1}>\alpha_{j,2}>\ldots>\alpha_{j,m}>\ldots>\alpha_{j,k}$.\smallskip
\noindent
Remark that for some decompositions of $w_0$, the map $\alpha_{j,m}\mapsto
\alpha_{j,m+1}$ corresponds to
the translation functor in the Auslander-Reiten quiver, [3].
\par\noindent
{\bf 1.6.} Let $(\,,\,)$ be the canonical pairing  between $\um$ and
$\up$,
see [4, 1.2],
 $\B$ be Lusztig's canonical base of $\um$, [{\it loc. cit.}],
and $\B^*\subset\up$ be the dual base, i.e. $(b^*,b')=\delta_{b,b'}$.
Let $u\mapsto \overline u$ be the antihomomorphism of
$\um$ such that $\overline F_i=F_i$ and $\overline q=q$.
Let $\tilde E_i$, $\tilde F_i$ : $\um\rightarrow\um$ be the Kashiwara
operators, [{\it loc. cit.}]. For $b\in\B$, $\tilde E_i(b)$, resp.
$\tilde F_i(b)$, equals some $b'$ in $\B\cup\{0\}$ modulo
$q^{-1}\Z[q^{-1}]\B$.
The rule $b\mapsto b'$ defines maps  $\tilde e_i$ and $\tilde f_i$ from
$\B$ to
$\B\cup\{0\}$.
For $b\in\B$, $1\leq i\leq n$, set $\varepsilon_i(b)=\hbox{Max}\{r,\,
\tilde e_i^r(b)\not=0\}$, and
$\E(b)=\sum_{i=1}^n\varepsilon_i(b)\varpi_i$.
The following is well
known, [24, Proposition 8.2], [18, 6.2.18-6.2.19] :\medskip\noindent
{\bf Theorem.} {\it For all $\lambda\in P^+$, we have :\par\noindent
\item {(i)} Via the isomorphism of Theorem 1.4 (iii),
$\ad\up K_{-2\lambda}$ is isomorphic to $\lql^*\otimes v_\lambda$
as a $\up$ module,\par\noindent
\item {(ii)} $\ad\up K_{-2\lambda}$ is generated as a space by $\{
K_{-2\lambda}b^*,\, b^*\in\B(\lambda)^*\}$, where $\B(\lambda)^*:=
\{b^*\in \B^*,\, \E(\overline b)\preceq\lambda\}$.}\qed\medskip\noindent
{\bf 1.7.} Recall that the enveloping algebra $U:=U(\g)$ is endowed
with a canonical filtration $\{U_k, k\in \N\}$ such that the associated
graded algebra is commutative, [14, 2.3]. Moreover, this filtration is
compatible with the $\ad U(\g)$-module structure of $U(\g)$. A
generalization of this filtration in the quantum case is the
Joseph-Letzter $P^+$-filtration, [21]. As in [27], this will
permit us to define the
generalized ($P^+$)-exponents of $\g$. First, let's present a few facts on
the
Joseph-Letzter filtration.\par\noindent
For $\lambda\in P^+$ set
$$F_q^\lambda:=\ad\uq.K_{-2\lambda},\hskip 1cm
F_{q,\lambda}:=\bigoplus_{\nu\leq\lambda}F_q^\nu.$$
Then, $F_{q,\lambda}$ is an $\ad\uq$-module.
By [18, 7.1.1], $\{F_{q,\lambda},\,\lambda\in P^+\}$ defines a
$P^+$-filtration
of $F(\uq)$. Moreover, let
$\hf^\lambda$, resp. $\hf(\uq)$, be the associated graded
space of $F_q^\lambda$, resp. $F_q(\uq)$. 
As  $\ad\uq$-modules, $\hf^\lambda$ and $F_q^\lambda$ are isomorphic.
By [20], the Joseph-Letzter's filtration extends to a filtration of
$\uq$ such that  $\hat U_q(\g)$ is the algebra
 with the canonical generators, quantum
Serre relations, weight relations, and
$[\hat E_i,\hat F_j]=-\delta_{ij}{\hat K_{-\alpha_i}\over q-q^{-1}}$.
Explicitly this filtration is ${\cal F}_\lambda(\uq)=\bigoplus_{\alpha\geq
-2\lambda}\up\otimes K_\alpha\otimes V_q(\nm)$
where $V_q(\nm)$ is the $\C(q)$-algebra generated by $K_{\alpha_j}F_j$,
$1\leq 
j\leq n$. In general, if $E$
is an element, resp. subspace, of $\uq$, $\hat E$ will be the associated
graded element, resp. space.
Clearly,
the algebra $\hat F(U_q(\g))\subset \hat U_q(\g)$ admits a Ringel
filtration as in 1.3 and the
associated graded algebra will be $\Gr\hf(\uq)$.
\par\noindent In $\hf(\uq)$ we have, [{\it loc. cit.}] :
$$\hf^\lambda \hf^{\lambda'}=\hf^{\lambda+\lambda'},\hskip 1cm \lambda,
\lambda'\in P^+.\leqno (1.7.1)$$
Hence, $\bigoplus_i\hf^{\varpi_i}$ generates the graded associated
algebra.\medskip\noindent
{\bf 1.8.}
For $\lambda\in P^+$, let $z_\lambda$ be the unique $\ad\uq$-invariant
element in $\ad\uq.K_{-2\lambda}$.  By Theorem 1.4, the center $Z_q$ of
$\uq$
is generated (as a
space) by the $z_\lambda$, $\lambda\in P^+.$ Moreover, $Z_q$ is generated
as
an algebra by the $z_i:=z_{\varpi_i}$, [18, 7.1.17].
\medskip\noindent
{\bf Definition 1.} Let $J_q$ be the ideal of $F(\uq)$ generated
by the $z_i-\varepsilon(z_i)$, $1\leq i\leq n$.
Let $\hat\hb$ be a complementary $\ad\uq$-module of $\hat J_q$ such that
$\hat\hb=\bigoplus_\lambda\hat\hb\cap \hat F_q^\lambda$. Let $\hb$ a
corresponding
sub-ad$\uq$-module of $F(\uq)$ such that $\hb=\oplus \hb\cap
F_q^\lambda$. For all $\lambda$ in $P^+$,  set $\hb^\lambda=\hb\cap
F_q^\lambda$, and let
 $\hb^\lambda(\mu)$, resp. $F_q^\lambda(\mu)$,
 be the isotypical component
of $\hb^\lambda$, resp. $F_q^\lambda$, of type $V_q(\mu)$.
\medskip\noindent
{\bf Definition 2.} Let $\mu\in P^+$. In the algebra $\C[P^+]$ generated
by
the $e^\lambda$,
$\lambda\in P^+$, we define the element $Q_\mu(e)=Q_\mu:=\sum_\lambda
[\hb^\lambda : V_q(\mu)]e^\lambda$. The dominant weights occurring  in
this
sum with multiplicities will be called $P^+$-exponents of
$\mu$.\smallskip\noindent
As in the classical case, we have, [18, Lemma 8.1.5] :
$$[\hb:V_q(\mu)]=\dim\,V_q(\mu)_0.\leqno (1.8.1)$$
Hence, by specializing $e$ on 1, we obtain $Q_\mu(1)=\dim\,V_q(\mu)_0$.
We now give an explicit expression of $Q_\mu$ in
terms of the canonical base. This is nothing but another formulation of a
result of Baumann [1, 3.4]. We give a proof for
completion.\medskip\noindent
{\bf Proposition.} {\it For $\mu$ in $P^+$,
let $B(\mu)$, resp. $B(\mu)_0$, be the crystal base,
resp. the set of
0-weight elements in the crystal base, of $V_q(\mu)$. Then,
$Q_\mu=\sum_{b\in B(\mu)_0}e^{\E(b)}$.\smallskip\noindent
Proof.} For $\lambda,\mu\in P^+$, consider the following subsets of
$B(\mu)_0$.
$$N_{\lambda,\mu}=\{b\in B(\mu)_0,\,\E(b)\preceq\lambda\},\hskip 2mm
 N_{\lambda,\mu}^0=\{b\in B(\mu)_0,\,\E(b)=\lambda\},\hskip 2mm
N_{\lambda,\mu}^+=N_{\lambda,\mu}\backslash N_{\lambda,\mu}^0,$$
with cardinals, respectively, $n_{\lambda,\mu}$,
$n_{\lambda,\mu}^0$ and $n_{\lambda,\mu}^+$.
 Now, let's work in the Joseph-Letzter
graded algebra. The theorem of separation of variables, [21], gives
$\hf=\hat\hb\otimes\hat Z_q$. Moreover, by semi-simplicity and
from the basic properties of the crystal base, [18, 6.3.18] :
$$[\hf^\lambda : V_q(\mu)]=[V_q(\lambda)^*\otimes V_q(\lambda):V_q(\mu)]
= [V_q(\lambda)\otimes V_q(\mu):V_q(\lambda)]=n_{\lambda,\mu}.
\leqno (1.8.2)$$ We claim that
$$\hf^\lambda(\mu)\simeq n_{\lambda,\mu}^0V_q(\mu)\oplus \hat J_q\cap
\hf^\lambda(\mu).
\leqno (1.8.3)$$
Indeed, we can prove this fact by induction on the $\vert\lambda\vert
:=\sum_i\lambda_i$, with $\lambda
=\sum_i\lambda_i\varpi_i$ in the following way.
Suppose this is true for all $\lambda'\in P^+$,
$\vert\lambda'\vert<\vert\lambda\vert$. Then,
we have [$\hat\hb(\mu)\cap\hf^{\lambda'} :
V_q(\mu)]=n_{\lambda',\mu}^0$, for such $\lambda'$. 
Note that $\hat J_q$ is the ideal generated by the
$\hat z_{\varpi_i}$, [18, 7.3.5]. By the separation of variables theorem
and
the induction hypothesis, this gives $$\sum_{\lambda'\prec\lambda}\,
\hat z_{\lambda-\E(b)}(\hat\hb\cap \hf^{\lambda'}(\mu))=\hat
J_q\cap\hf^\lambda(\mu)
.$$
By  separation of variables and the induction hypothesis, the left hand
term is isomorphic to
$\sum_{\lambda'\prec\lambda}n_{\lambda',\mu}^0V_q(\mu)
=n_{\lambda,\mu}^+V_q(\mu)$. By (1.8.2), this gives (1.8.3) and the
proposition is proved.\qed\medskip\noindent
{\bf Remark.} When $\g$ is of type $A_n$, then the function $Q_\mu$ is
the multivariable Kostka-Foulkes polynomial ${\K}_{\mu,0}$
defined in [29, 6.2] by
Lascoux-Leclerc and Thibon. In fact,
the previous proposition materializes their definition in the following
sense
:  ${\K}_\mu$ as defined in [{\it Loc. cit.}] is the Hilbert polynomial
of the isotypical component of $\hb$ of type $V_q(\mu)$ in the
Joseph-Letzter
$P^+$-filtration. \medskip\noindent
{\bf 1.9.} We present a result which can be seen as a quantum version
of Dixmier's antihomomorphism [14, 8.4.1].\par\noindent
We consider the Joseph-Letzter associated graded algebra.
Let  $\hat\pi$ be the natural projection from $\hf(\uq)\simeq
\oplus_\lambda\lql^*\otimes\lql$ onto $\oplus \ad\up\,\hat
K_{-2\lambda}\simeq
\oplus\lql^*\otimes v_\lambda$. Then, by [22, Lemma
7.2],\smallskip\noindent
{\bf Lemma.} {\it The restriction of $\hat\pi$ to $\hf(\uq)^{\up}$ is
an injective algebra $q$-anti\-ho\-mo\-mor\-phism. To be more precise, let
$a_\mu
\in\hf(\uq)^{\up}$ ,
$b_\nu\in(\hf^{\lambda})^{\up}$ of weights, respectively, $\mu$ and $\nu$.
Then, $\hat\pi(a_\mu
b_\nu)=q^{(\mu,2\nu-\lambda)}\hat\pi(b_\nu)\hat\pi(a_\mu)$.
Moreover $\hat\pi(\hat z_\lambda)=\hat K_{-2\lambda}$.}\qed
\medskip\noindent
For $\lambda$ in $P^+$, let $\hat\pi_\lambda$ be the $\lambda$ component
of 
$\hat\pi$.
By say [8, Proposition 3.4], it is clear that $\hat\pi$ is the restriction
of
the natural projection on $\hat\up\otimes \hat K_{-2\lambda}$.
Hence, $\hat\pi$ is compatible with the Ringel
filtration.\medskip\noindent
{\bf 2. $A$-form and specialization.}\bigskip\noindent
{\bf 2.1.} The following results can be deduced from [12, 1.5].
Let $\A$ be the local algebra $\C[q]_{(q-1)}$. We define the
$\A$-form $\ua$ generated by $E_i$, $F_i$,
${K_{\varpi_i}-K_{-\varpi_i}\over q-q^{-1}}$, $K_{\lambda}$, where
$\lambda\in P$, $1\leq i\leq n$. Remark that $\ua$ is a
sub-$\ad\,\ua$-module
of $\uq$. Moreover, $\ua$ is $\A$-free and $\uq=\C(q)\otimes_{A}\ua$.
\medskip\noindent
{\bf Claim.} {\it Let $\C[P/2P]$ be the algebra of the group $P/2P$,
generated
as a space by the elements $e^{\overline\lambda}$, $\overline\lambda
\in P/2P$.
There exists a natural isomorphism between $\ua/(q-1)\ua$ and
the central extension $\C[P/2P]
\otimes U(\g)$ such that the canonical surjection
$\phi\,:\,\ua\rightarrow\,
\ua/(q-1)\ua\simeq 
\C[P/2P]
\otimes U(\g)$ sends $K_\lambda$ to $e^{\overline\lambda}$,
$\lambda\in P$.}\qed\smallskip\noindent
 Let $\uap$, resp. $\uao$, resp. $\uam$, be the positive,
resp. Cartan, resp. negative, part of $\ua$. As in 1.2, the triangular
decomposition holds for $\A$-forms.
\par\noindent We define $F(\ua)=\ua\cap F(\uq)$ and
$\za:=Z\cap\ua$. From the definitions, the algebra $\za$ is
the center of $\ua$.
Let $J_q$ be the ideal of $F(\uq)$ generated
by the $z_i-\varepsilon(z_i)$, $1\leq i\leq n$ and set $\ja:=J_q\cap
F(\ua)$.
Then, $\ja$ is a
$\ad\ua$-module. Moreover, the quotient $F(\ua)/\ja\subset
F(\uq)/J$
has no $\A$-torsion. 
\smallskip\noindent
{\bf Proposition.} {\it The morphism $\phi$ sends 
\item{(i)} $F(\ua)$ onto $U(\g)$, realized as the subalgebra $1\otimes U(\g)$ in $\C[P/2P]
\otimes U(\g)$, 
with kernel $(q-1)F(\ua)$.\par\noindent
\item{(ii)}  $\ja$ onto the minimal primitive ideal
$U(\g)\Ker\varepsilon_1$, where
$\varepsilon_1$ is the augmentation restricted to $Z(\g)$. 
\smallskip\noindent
Proof.} By Theorem 1.4 (ii), $F(\ua)$ contains
$K_{-2\varpi_i}E_{\alpha_i}$,
for all $i$. Hence, by the Claim, $\phi(F(\ua))$ contains $U(\n)$. In the same
way,
it contains $U(\nm)$, and both algebras generate $U(\g)$. Hence,
$U(\g)\subset
\phi(F(\ua))$. Now, by [20, Theorem 6.4],
$F(\ua)\subset \uam\otimes \hat U^0\otimes
\uap$, with $$\hat U^0=\C(q)[K_{2\varpi_i},K_{-2\varpi_i}]\cap
\A[K_{\varpi_i},K_{-\varpi_i},
{K_{\varpi_i}-K_{-\varpi_i}\over
q-q^{-1}}]=\A[K_{2\varpi_i},K_{-2\varpi_i},
{1-K_{-2\varpi_i}\over q-q^{-1}}].\leqno (2.1.1)$$
Hence, $\phi(F(\ua))\subset U(\g)$. This gives (i).\par\noindent
Let's prove that $\phi(\ja)=U(\g)\Ker\varepsilon_1$.\par\noindent
$U(\g)\Ker\varepsilon_1\subset\phi(\ja)$ : Let $J_{\A}'=
\fa(\Ker\varepsilon\cap\za)$. Then, $\ja\subset J_{\A}'$ implies
$U(\g)\Ker\varepsilon_1
=\phi(J_{\A}')\subset\phi(\ja)$.\par\noindent
$\phi(\ja)\subset U(\g)\Ker\varepsilon_1$ : Let $M_q(0)$ be the quantum Verma
module with highest weight 0, as defined in [18, 3.4.9]. Then, $M_q(0)$ is
a cyclic
$\um$-free module generated by a highest weight vector $v_0$ of weight 0.
Let $M_{\A}(0)$ be the sub-$U_A$-module of $M_q(0)$ generated by $v_0$. Then,
$M_{\A}(0)/(q-1)M_{\A}(0)$ has a natural structure of $U(\g)$-module by (i).
It is easily seen that this is a cyclic $U(\g)$-module with highest weight 0 and
with same formula character than the classical Verma module $M(0)$. Hence,
$M_{\A}(0)/(q-1)M_{\A}(0)\simeq M(0)$ as $U(\g)$-modules. Now,
$\ja$ annihilates $M_{\A}(0)\subset M_q(0)$. Hence,
$\phi(\ja)\simeq\ja/(q-1)\ja$ annihilates $M(0)$. By Duflo's theorem [14,
8.4.3],
we have the desired inclusion. Hence, (ii) holds.
\qed
\medskip\noindent
{\bf 2.2.}  Let $\up^+$, resp. $\um^+$, be
the augmentation ideal of $\up$, resp. $\um$.
Let $\varphi$ be the Harish-Chandra map of $\uq$, i.e. the projection $\uq
\rightarrow \uo$ with kernel $\uq\up^++\um^+\uq$. Recall
that $\uo$ is isomorphic to $\C(q)[P]$ as a $W$-module. 
Recall the so-called twisted action of $W$ on $\uo$, see [12, §2].
Moreover, by [12, 2.2] :\medskip\noindent
{\bf Proposition.} {\it Let $\tilde U_{\A}^0$
be the intersection of $\ua$ with the $\C(q)$-algebra generated by the
$K_{2\lambda}$, $\lambda\in P$.  Then, $\tilde U_{\A}^0$ is a
sub-$W$-module
of $U^0$ and the Harish-Chandra homomorphism
$\varphi$ restricts into an algebra isomorphism $Z_{\A}\simeq
(\tilde U_{\A}^0)^W$.}\qed
\medskip\noindent
Remark that the inverse isomorphism can be given explicitly by the Tolstoi
projector
 introduced in [26]. This is implicit in the proof of [12, Proposition
2.2].
The proposition above is the quantum integral analogue of the following
classical result, see
[14, 7.4] : the Harish-Chandra map etablishes an algebra isomorphism
between the center $Z(\g)$ of $U(\g)$ and the
$W$-invariant subalgebra of the Cartan part $U(\h)$ of
$U(\g)$.\medskip\noindent
{\bf Theorem.} {\it The $A$-algebra $Z_A$ is a polynomial algebra over $A$. 
It specializes for
$q=1$, onto the center $Z(\g)$ of the enveloping algebra
$U(\g)$, realized as the subalgebra $1\otimes U(\g)$ in $\C[P/2P]
\otimes U(\g)$.\smallskip\noindent
Proof.} From Proposition 2.1, $\za$ specializes in the subalgebra  $1\otimes U(\g)$.
Clearly, specialization sends an element of the center
of $\ua$ to an element of $Z(\g)$. \par\noindent
Moreover, $\tilde U_{\A}^0$ specializes onto $U(\h)$ and specialization
commutes
with $W$-action.
As $W$ is finite, the same is true for the $W$-invariants. The last assertion of
the theorem
follows because specia\-li\-za\-tion commutes with the Harish-Chandra
homomorphisms.\par\noindent
Let's prove that $\za$ is polynomial over $A$.\par\noindent
Let's identify the center $Z$ of $U(\g)$ with $S(\h)^W$ and the 
symmetric algebra $S(\h)$ with the algebra of regular functions on the
lattice of radical weights $Q$. In the same way, let's identify $Z_q$ with 
$(\uqo)^W$, and $\uqo$ with the 
$\C(q)$-algebra generated by the functions $K_{\varpi_i}(\_)=q^{(\varpi_i,\_)}$
 on $Q$.
This enables us to give formal expansions at $q=1$ from elements of
 $Z_q$, as in [8, 7.1] and 
 [JL2, 6.13], i.e. we can embed $Z_q$ in $S(\h)^W[[q-1]]\simeq Z[[q-1]]$. 
 Let's consider the augmentation ideal $Z^+$ 
of $Z$. Fix a set of homogeneous elements $(C_{m_j})$, $1\leq j\leq n$,
of degree $0<m_1\leq \ldots\leq m_n$, which is a base of $Z^+$ modulo $(Z^+)^2$. 
Let $z_i$ be the quantum central element as in 1.8 and 
$z'_i=z_i-\varepsilon(z_i)$.
From the formula in the proof of [JL2, Lemme 6.15], we have the claim \medskip\noindent
{\bf Claim.}  {\it The coefficient of $(q-1)^m$ of the image of $z'_i$ in
$Z[[q-1]]$ 
is an homogeneous element of degree $m$ in $Z$.}\qed\medskip\noindent 
Now, the image of $z'_i$ in $Z^+/(Z^+)^2[[q-1]]$ is $\sum 
a_{ij} C_{m_j}(q-1)^{m_j}$, $a_{ij}\in\C$, $1\leq i,j\leq n$. 
Moreover the products $z'_iz'_j$ is zero in this quotient.
Recall that the $z'_i$ generate the $\C(q)$-algebra  $Z_q$ and recall that, by the previous arguments,
$Z_A$ specializes onto $Z$ at $q=1$. This implies that there exists 
a polynomial combination of the  $z'_i$ (and of $(q-1)^{-1}$) which specializes 
on  $C_{m_j}$ at $q=1$. This combination is linear modulo $(Z^+)^2$
and this implies that the matrix $(a_{ij})$ is invertible. Hence, we can change the 
index of the $z_i$ in order to obtain non zero minors  
$(a_{ij})_{1\leq i,j\leq k}$, $1\leq k\leq n$.  Then, by induction on $k$, we obtain 
elements 
$C_{m_k}^q$ in $Z_A$ which are linear combination of $z_k$ and products of
 $z_i$ for $i<k$ and which specializes on $C_{m_k}$. 
 Now, we assert that $Z_A$ is the $A$ polynomial algebra generated by the
  $C_{m_k}^q$. Indeed, as in the proof of Lemma A1, 
  this is a direct consequence of the following facts :
  \medskip\noindent
1) These elements specialize on a polynomial base of $Z$.\smallskip\noindent
2) They generate $Z_q$ as a $\C(q)$-algebra. 
\qed
\medskip\noindent
{\bf 2.3.} 
\def\fa{F(\ua)}\def\far{F_{\A,\lambda}}\def\faru{F_{\A,\mu}}
\def\fal{F_{\A}^\lambda}
\def\fua{F(\ua)}
 Set $\far=\ua\cap F_{q,\lambda}$.
\medskip\noindent 
{\bf Proposition.} {\it We have the following\par\noindent
\item{(i)} The $\A$-module $\fa$ is free,\par\noindent
\item{(ii)} there exists
an ad$\ua$-modules graduation $(\fal)_{\lambda\in P^+}$ 
of $\fua$ such that
$$\fua=\oplus_{\lambda\in P^+ }\fal,\hskip 5 mm 
\oplus_{\lambda\leq \mu}\fal=\ua\cap
F_{q,\mu}.$$ \smallskip\noindent
Proof.} The module $\ua$ is $\A$-free. In particular,
 $\far$ has no $\A$-torsion, so it is a free $\A$-module.
Fix $\lambda\leq\mu$ in $P^+$.
Then, $\faru/\far$ embeds in $F_{q,\mu}/F_{q,\lambda}$, so it has no
$\A$-torsion.
Hence, the $\A$-module $\far$ is a direct summand in $\faru$. Remark that
$\fa=\bigcup_\mu\faru$ to obtain (i) and (ii).\qed\medskip\noindent
{\bf Remark.} 
$F(\ua)$ strictly contains
$\oplus_\lambda (\ua\cap F_{q}^{\lambda})$. For example, set $\g=\sl_2$ and let
$z_\varpi\in F_{q}^{\varpi}$ be the Casimir element corresponding to
the quantum trace of the fundamental module $V_q(\varpi)$. Then,
[8, 5.3], $\phi(z_\varpi)=2$, ${z_\varpi-2\over q-q^{-1}}$ belongs to
$F(\ua)$ but not to $\oplus_\lambda \ua\cap F_{q}^{\lambda}$.
\medskip\noindent
{\bf 2.4.} We prove in this section that the separation theorem at the
quantum level, i.e. on $\C(q)$, and on the classical level implies 
the separation theorem for the $\A$-forms. The proof relies on general
properties of torsion-free modules over $\A$-polynomial algebras. These properties
are given in the Appendix. \smallskip\noindent
{\bf Theorem.} {\it The separation theorem holds for the $\A$-forms, i.e.
via multiplication
$$\za\otimes\ha\simeq\fa,$$  for an ad$\ua$-module
$\ha$. Moreover, $\ha$ can be obtained as $\ua\cap\hq$, where $\hq$ is as 
in [1.8, Definition 1].
\smallskip\noindent
Proof.} Fix $\hq$ as in Definition 1 and set $\ha=\ua\cap\hq$.  
It is sufficient to prove the theorem on the isotypical component,
i.e
$F_A(\lambda)
=Z_A\otimes H_A(\lambda)$, $\lambda\in P^+$. Now, from Theorem 2.2, 
the center $Z_A$ is a polynomial ring over $A$. 
We are in the framework of 
[Appendix, Proposition A4] 
with 
 $T=q-1$ and $R=Z_A$. From the qunatum and the classical 
 theorem of separation of variables, [21] and [27], the
hypothesis of Proposition A4 are verified with $r=$dim$V(\lambda)$dim$V(\lambda)_0$.
Hence, $F_A(\lambda)$ is a $Z_A$-free module with rank $r$. 
\smallskip\noindent
Fix a base of this module.
The theorem is
a consequence of the following facts :\smallskip\noindent
(i) $H_A(\lambda)$ generates the $Z_q$ module $F_q(\lambda)$.\par\noindent
(ii) $H_A(\lambda)$  is a $A$-direct summand in 
 $F_A(\lambda)$.\smallskip\noindent
Indeed, fix a $A$-base of  $H_A(\lambda)$  which is compatible with the
graduation
of  $F_A(\lambda)$. Then, this base generates the $Z_q$ module 
$F_q(\lambda)$ by [18, 7.3.8]. This implies (i).\par\noindent
Moreover, by contruction $\ha(\lambda)$ is a free $A$-module with finite rank
and \par\noindent $\fa(\lambda)/\ha(\lambda)$ has no $(q-1)$-torsion. So, (ii) holds.
\par\noindent
Now, fix a $A$-base of $\ha(\lambda)$ and
let $d$ be the  determinant of the matrix of this base in
the fixed base of the $Z_A$-module $F_A(\lambda)$. Then, (i) asserts that  
$(q-1)^md\in A^*$ for some positive integer $m$. By (ii), we obtain $m=0$.
This ends the proof. \qed\medskip\noindent
\medskip\noindent
Let $\phi$ be the morphism of specialization as in 2.1. We have
:\medskip\noindent
{\bf Corollary.} {\it The morphism $\phi$ sends \par\noindent
\item{(i)} $\ha$ onto a complementary subspace  of $U(\g)\Ker\varepsilon_1$ in
$U(\g)$,
\item{(ii)} The $\mu$ isotypical component of $\ha$
onto the $\mu$ isotypical component of some module of harmonic elements of
$U(\g)$ (in the sense of 0.4).
\smallskip\noindent
Proof.} Remark that $\za=(\za\cap\Ker\varepsilon)\oplus A$ as $A$-modules. 
By the previous theorem, this implies that $\fa=\ja\oplus\ha$. Hence,
Proposition 2.1 gives the corollary.\qed\medskip\noindent

{\bf 2.5.} From Theorem 2.4, we can construct a graduation 
$(\fal)_{\lambda\in P^+}$ such that $\ha=\bigoplus\ha\cap\fal$
and $\za=\bigoplus\za\cap\fal$. Let $Z_A^\lambda$ and
$H_A^\lambda$ be the graded component of $Z_A$ and $H_A$. We have 
$$\bigoplus_{\mu'\preceq\mu\leq\lambda}\za^{\mu-\mu'}\otimes
\ha^{\mu'}=F_{\A,\lambda}.\leqno (2.5.1)$$
For all $\lambda\in P^+$,
set $F_\lambda=F_\lambda(U(\g))=\phi(\far)\in U(\g)$. This is an
ad$U(\g)$-submodule of $U(\g)$. In general, we omit the index $A$ to note the image by $\phi$,
i.e. we note $E:=\phi(E_A)$. We want to prove
that $F_\lambda$ is a filtration of $U(\g)$ which is analogue to the
Joseph-Letzter's
filtration in the quantum case. Indeed we have :\smallskip\noindent
{\bf Proposition.} {\it The algebra $U(\g)$ is filtered by the
ad$U(\g)$-modules $(F_\lambda)_{\lambda\in P^+}$.
Moreover, \smallskip\noindent
\item{(i)} $F_\lambda\simeq\bigoplus_
{\mu\leq\lambda} F^\mu$\par\noindent
\item{(ii)}
$F_\mu\simeq V(\mu)^*\otimes V(\mu)$\par\noindent
\item{(iii)}
$\bigoplus_{\mu'\preceq\mu\leq\lambda}Z^{\mu-\mu'}\otimes
H^{\mu'}=F_\lambda.$\smallskip\noindent
Proof.} The ad$\ua$ modules $(\far)_{\lambda\in P^+}$ give a filtration of
$\fa$. Hence, by specialization and by Proposition 2.3, we have the first
assertion. (i) and (ii) are obtained by specialization, using the fact that
tha $A$-spaces $\far$ and $F_A^\lambda$ are $A$-summands in $\fa$.\par\noindent
(iii) is obtained from (2.5.1).\qed
\medskip\noindent
{\bf 3. Gr\"obner bases and harmonic elements.}\bigskip\noindent
{\bf 3.1.}  Recall the definition of roots packages, see 1.5. In the
graded
algebra $\Gr\uq$,
the generators $\Gr\,z_i$ of $\Gr\,Z_q$ are given by :\medskip\noindent
{\bf Proposition.} {\it For $1\leq i\leq n$, we have, up to a
multiplicative
scalar
$$\Gr\,z_i=(\prod_{\beta\in\Delta_i^+}\Gr\,E_\beta
\Gr\,F_\beta)\Gr\,K_{-2\varpi_i}.$$
The $\Gr\,z_i$, $1\leq i\leq n$, generate $\Gr\,Z_q$.
\smallskip\noindent
Proof.} Recall that $z_i$ corresponds to the quantum trace on
$V_q(\varpi_i)$, [8, 5.1].
The formula is a direct consequence of [8, Proposition 3.4], and
[11, Lemma 4.1 (iii)]. From 1.5, $\{\Delta_i^+,\, i\}$ is a partition
of $\rp$, so the $\Gr\,\prod_iz_i^{a_i}$, $a_i\in\N$ are linearly
independant.
Hence, they generate $\Gr\,Z_q$ as a space. This proves the proposition.
\qed\medskip\noindent
{\bf 3.2.} \par\noindent
{\bf Definition.} Let $\M$ be the set of monomials $M$ in
$\uq$ such that \smallskip\noindent
\item{(i)} $M=\prod_{j,m}E_{\alpha_{j,m}}^{a_{j,m}}
K_\lambda F_{\alpha_{j,m}}^{b_{j,m}}$, $0\leq a_{j,m}, b_{j,m}$,
$\lambda\in P$,\par\noindent
\item{(ii)} For all $j$, $\prod_ma_{j,m}b_{j,m}=0$.\smallskip\noindent
Then, the following separation theorem holds.\medskip\noindent
{\bf Theorem.} {\it Let $J_q^0$ be the ideal of $\uq$ generated by the
$z_i$, $1\leq i\leq n$. Fix a reduced decomposition of $w_0$.  Then
:\par\noindent
\item{(i)} The set $\{z_i, \, 1\leq i\leq n\}$ is
a Gr\"obner base of $J_q^0$ for the lexicographic ordering of the
PBW-base, i.e. $u\in J_q^0=<z_i>\Rightarrow
\Gr\,u\in<\Gr\,z_i>$.\par\noindent
\item{(ii)}  The algebra $\Gr\,\uq$
is free over $\Gr\,Z_q$ with base $\Gr\,\M$.\par\noindent
\item{(iii)} The algebra $\uq$ is free over  $Z_q$ with base $\M$.
\smallskip\noindent
Proof.} First remark that, by construction and by Proposition 3.1,
the image of $\Gr\,\M$ in
$\Gr\,\uq/<\Gr\,z_i, 1\leq i\leq n>$ is a base of this
quotient.\par\noindent
Let $\Gr\,z_j$, resp. $\Gr\,m_j$, be monomials in $\Gr\,Z_q$, resp.
$\Gr\,\M$,
$j=1,2$. Then,
$\{\Delta_i^+, i\}$ being a partition of $\rp$, we obtain by Proposition
3.1 :
$$\Gr\,z_1 \Gr\,m_1=\Gr\,z_2 \Gr\,m_2 \Rightarrow
\Gr\,z_1 =\Gr\,z_2 \;\hbox{ and }\;  \Gr\,m_1=\Gr\,m_2.\leqno (3.2.1)$$
Now, considering $\Gr\,\uq$ as a $\N^N\times P\times \N^N$-graded algebra,
see 1.3, we remark that each graded component of $\Gr\,\uq$ has dimension 0 or 
1. Hence,
 (3.2.1) implies that $\Gr\,\uq$
is free over $\Gr\,Z_q$, with base $\Gr\,\M$. This implies (ii) and then
(iii).\par\noindent
Let $Z_q^+$ be the ideal of $Z_q$ generated by the $z_i$. Then,
by (iii), each element $u$ of $J_q^0$ can be decomposed in a unique way
into
$u=\sum_kP_km_k$, with $P_k\in Z_q^+$ and $m_k\in\M$. This implies
(i) by (3.2.1).\qed \smallskip\noindent
{\bf Remark 1.} It is possible to generalize the theory of
Gr\"obner bases in some non-com\-mu\-ta\-ti\-ve contexts, [23].
(i) can also be proved in this context by using the first Buchsberger
criterion [2]. Indeed,
Proposition 3.1 asserts that the initial monomials of the $z_i$ are
relatively
 prime.\smallskip\noindent
{\bf Remark 2.} Unfortunately, the specialization of the center
is not compatible with the associated graded algebra. Indeed, the initial
terms
of the generators $z_i$ may vanish at $q=1$. The base given in Theorem 3.2
(iii) has no analogue in the classical case.\bigskip\noindent
{\bf 3.3.} This section will be of no use in the sequel.
The results of 1.3 generalize when $q=\varepsilon$ is a root of one.
We define the quantum algebra $U_\varepsilon(\g)$ as in [12, 1.5].
By the description of the center $Z_\varepsilon$ of quantum groups at
roots
of one,
see [13], we obtain :
\medskip\noindent
{\bf Theorem.} {\it When $q=\varepsilon$ is a  root of one, the algebra
$\Gr U_\varepsilon(\g)$ is free over the graded algebra
$Gr\,Z_\varepsilon$.
This implies that $U_\varepsilon(\g)$ is free over
its center.}\qed\smallskip\noindent
For example, when $\varepsilon$ is a $l$th root of one, with $l$ odd,
then a base
of $U_\varepsilon(\g)$ over $Z_\varepsilon$ is given by $\M^l:=\{M\in\M,
M=\prod_{i,j,m}E_{\alpha_{j,m}}^{a_{j,m}}
K_{\varpi_i}^{c_i} F_{\alpha_{j,m}}^{b_{j,m}}, 0\leq a_{j,m},
\,b_{j,m},\,c_i<l\}$.\medskip\noindent
{\bf 3.4.} Let $\pi$ as in 1.9. The next proposition gives
 a base of $\pi(\hb^{\up})$ in terms of
the dual of the canonical base of $\um$, (see 1.6). When $u$ is a weight
vector in a
$\uq$-module, then we note $\Omega(u)$ the weight of $u$.
\medskip\noindent
{\bf Proposition.} {\it For $\lambda$ in $P^+$, set
$$\B_H:=
\{b, \,b\in \B,\,\E(b)\preceq -\Omega(b)\},\hskip 5 mm
  \B_H^\lambda:=
\{b, \,b\in \B_H, \E(\overline b)=\lambda\}.$$
Then, a base of
$\pi(\hb^{\up})$ is given by $(K_{-2\lambda}b^*_\lambda,\,\lambda\in P^+,
\,b_\lambda\in \B_H^\lambda)$ for some harmonic space
$\hb$.\smallskip\noindent
Proof.} Fix $\mu=\sum \mu_i\varpi_i\in P^+$.
Recall the notations of Proposition 1.8. we have
 $$\dim\,(\hb^{\lambda})(\mu)^{\up}=
[\hb^\lambda(\mu):V_q(\mu)]=n_{\lambda,\mu}^0.\leqno (3.4.1)$$
Via the natural embedding of $B(\mu)$ in the canonical base,
we have by [24, Proposition 8.2]
 $N_{\lambda,\mu}^0=\{b,\,b\in \B,\,\Omega(b)=-\mu, \E(\overline
b)\preceq\mu,\,
\E(b)=\lambda\}$. In the same way,
dim$(\hf^{\lambda})^{\up}(\mu)=n_{\lambda,\mu}$ and
$N_{\lambda,\mu}=\{b,\,b\in \B,\,\Omega(b)=-\mu, \,\E(\overline
b)\preceq\mu,
\E(b)\preceq\lambda\}$.\smallskip\noindent
{\bf Claim.} {\it A base of
 $\pi(F_q^\lambda(\mu)^{\up})$ is given by $(K_{-2\lambda}b^*)$,
where $b$ runs over the set $\overline N_{\lambda,\mu}$, i.e. the image of
$N_{\lambda,\mu}$ by the involution ${\,}^{\overline{\hbox{\hskip 2mm}}}$
described in 1.6.
\smallskip\noindent
Proof of the claim.} Fix $\lambda$, $\mu$, $\xi$ in $P^+$,
$\lambda=\sum\lambda_i\varpi_i$,
$\mu=\sum\mu_i\varpi_i$. By say, [18, 6.3.20]
the image of the highest weight
space of weight $\mu$ in $V_q(\xi)\otimes\lql$ by the natural projection
on $V_q(\xi)\otimes v_\lambda$ is $(V_q(\xi))_{\mu-\lambda}^\lambda\otimes
v_\lambda$, where
$$(V_q(\xi))_{\mu-\lambda}^\lambda:=\{v\in V_q(\xi),
\,\Omega(v)=\mu-\lambda,\,
F_i^{\mu_i+1}v=0, \forall i\}\hskip 3cm$$
$$\hskip 5cm =\{v\in V_q(\xi), \,\Omega(v)=\mu-\lambda,\,
E_i^{\lambda_i+1}v=0, \forall i\}.\leqno (3.4.2)$$
We have the succession of space isomorphisms
$(\lql^*)_{\mu-\lambda}^\lambda=<b^*\in \B(\lambda)^*, \,\Omega(b)=-\mu,\,
\E(b)\preceq\mu>=<b^*\in \B^*, \,\Omega(b)=-\mu,\,
\E(\overline b)\preceq\lambda,\,\E(b)\preceq\mu>=\overline
N_{\lambda,\mu}$.
Indeed, the first equality comes from [30, Theorem 4.4 (c), 4.3] and
(3.4.2),
the second one is deduced from [24, Proposition 8.2].
The claim results from Theorem 1.4 (ii) and Theorem 1.6.
\qed
\smallskip\noindent
{\it end of proof of the proposition.} Suppose $\nu$ in $P^+$,
$\nu\preceq\lambda$.
 In the graded algebra  $\hf(\uq)$, $<
K_{-2\lambda}b^*,\, b^*\in\B(\nu)_\mu^*>\subset \hat\pi(\hat J_q)$ by the
claim
and Lemma 1.9.
 Then, the proposition results from the claim by
(3.4.1).\qed\medskip\noindent
{\bf Remark 1.} What is behind the proof is a natural twist $\tau$ :
$B(\lambda)\otimes B(\mu)\rightarrow B(\mu)\otimes B(\lambda)$
defined to be the unique crystal map such that
 the image of the $\tilde e_i$-invariant elements, $1\leq i\leq n$, is
given by
$\tau(v_\lambda\otimes b v_\mu)=v_\mu\otimes \overline b v_\lambda $,
for $\E(b)\preceq\lambda$, $\E(\overline b)\preceq\mu$. It would be
interesting
to understand the braiding of this twist on $B(\varpi)^{\otimes k}$,
where $\varpi$ is a basic weight.\medskip\noindent
{\bf Remark 2.} We can prove that the generators of $\hb^{\up}$ can
be calculated from the $K_{-2\lambda}b^*_\lambda,\,\lambda\in P^+,
\,b_\lambda\in \B_H^\lambda$, by applying the Tolstoi projector, see [26].
\medskip\noindent
{\bf 3.5.} We suppose in this section that $\g$ is of type A-D-E
and that the decomposition of $w_0$ corresponds to an orientation of the
Coxeter graph of $\g$. In this case, from Ringel's Hall
algebra approach, [33],
the canonical base can be labelled by the Kostant partitions.
We present a few results
which can be found in [31], [32, par 2 and 3].
Recall that $\{E_\psi,\;\psi\in \N^N\}$ is a base of $\up$. There exists
 a parametrization $\B=\{b_\psi, \psi\in \N^N\}$ of the  canonical base
such that :
$$\Gr\,b_\psi^*=
\Gr\,E_\psi, \,\,\hbox{ up to a multiplicative scalar in }
\Z[q],
\leqno (3.5.1)$$

$$\Gr\,b_\psi^*b_\phi^*=\Gr\,q^{[\psi,\phi]}b_{\psi+\phi}^*,\,\,
\hbox{for a power $[\psi,\phi]$ of } q,
\leqno (3.5.2)$$

$$\E(b_{\psi+\phi})\preceq \E(b_{\psi})+\E(b_{\phi}).\leqno (3.5.3)$$
The following definition flows from Proposition 3.4 :\medskip\noindent
{\bf Definition.} The harmonic cone is defined by $\co:=
\{\psi,\,b_\psi\in \B_H\}\subset\N^N$.\medskip\noindent
For $\psi$ in $\N^N$, set $\Omega(\psi)=\sum_{l=1}^N \psi_l\beta_l$ be the
weight of
$E^{\psi}$ and $\E(\psi)=\E(\overline b_\psi)$.
Let $<\,,\,>$ be the bilinear form defined by the $\Delta^+\times
\Delta^+$
matrix $M=(m_{\alpha,\beta})$, with $$m_{\alpha,\beta}=
\cases{0 & if $\beta\leq\alpha$,\cr
(\alpha,\beta) & if not}$$
{\bf Theorem.} {\it The algebra $\C[\co]$ of the semigroup $\co$ is
normal.
Via the anti-ho\-mo\-mor\-phism
$\Gr\,\hat\pi$, see 1.9, the algebra $\Gr\,(\hf(\uq)/\hat J)^{\up}$
 is isomorphic to the following algebra : a base is
$a_\psi$, $\psi\in \co$ with the multiplication given
by $$a_\psi a_\phi=
\cases{0 & if $\E(\psi+\phi)\prec\E(\psi)+\E(\phi)$,\cr
q^{<\psi,\phi>-2(\Omega(\psi),\E(\phi))}a_{\psi+\phi} & if
$\E(\psi+\phi)=\E(\psi)+\E(\phi)$.\cr}$$
\smallskip\noindent
Proof.} The cone $\co$ is
$\{\psi\in\N^N,\;\E(b_\psi)\preceq\Omega(\psi)\}$.
$\Omega$ is clearly additive on $\N^N$ and so, $\co$ is a semigroup by
(3.5.3).
For all $\psi$ in $\N^N$, set
 $a_\psi:=\Gr\,\hat\pi(\Gr\,\hat K_{-2\E(\psi)}\hat E^\psi)$. The set
$\{a_\psi,
\,\psi\in\co\}$ is a base of $\Gr\,\hat\pi (\Gr\,(\hf(\uq)/\hat J)^{\up})$
by Proposition 3.4 and the previous formulas. The multiplication rule
is given by Lemma 1.9 and (1.3.2). \par\noindent
By [5, Theorem 3.8], we obtain that $\E(k\psi)=k\E(\psi)$, $\psi\in\N^N$.
 . This gives
$k\psi\in\co$, $k\in\N^*$, $\psi\in\N^N$ $\Rightarrow$ $\psi\in\co$. 
Hance, $\C[\co]$ is normal.\qed
\medskip\noindent
{\bf 4. The A$_n$ case.}\bigskip\noindent
{\bf 4.1.} Let $\g$ of type A$_n$. Set $F_q=F_q^{\varpi_1}$. For all
$\lambda
=\sum\lambda_i\varpi_i\in P$, we define $h(\lambda)=\sum_i i\lambda_i$.
Define the total order $\leq$ on $P^+$ such that \par\noindent
1) $h(\lambda)<h(\mu)\Rightarrow
\lambda\leq \mu$ and \par\noindent
2) Fix $\lambda=\sum_i\lambda_i\varpi_i$, $\mu=\sum_i\mu_i\varpi_i$, with
$h(\lambda)=h(\mu)$. Then, $\lambda\leq \mu$ iff $(\lambda_1,\ldots,
\lambda_n)$ is lower than  $(\mu_1,\ldots,
\mu_n)$ for the reverse lexicographical order of $\N^n$.\smallskip\noindent
{\bf Example.} For $\g=\sl_4$, $0\leq \varpi_1\leq \varpi_2\leq 2\varpi_1
\leq \varpi_3\leq \varpi_1+\varpi_2\leq 3\varpi_1\leq \varpi_1+\varpi_3
\leq 2\varpi_2\leq 4\varpi_1\leq \ldots$\smallskip\noindent
It is easily seen that this ordering verifies the hypothesis of 1.1.\par\noindent
We have :\medskip\noindent
{\bf Lemma.} {\it For all $k$ in $\N$, $\sum_{0\leq m\leq k} F_q^m=
\bigoplus_{h(\lambda)\leq k} F_q^\lambda$.\smallskip\noindent
Proof.} It is well known that the minimal integer $k$ such that
$\lql\subset V_q(\varpi_1)^{\otimes k}$ is $h(\lambda)$. Hence the lemma
results from
[19, Corollary 3.10].\qed\medskip\noindent
{\bf 4.2.} This section is devoted to the specialization of the center of
$U_q(\sl_{n+1})$. Recall that the center $Z_q$ of $U_q(\sl_{n+1})$ is
generated as a space
by the $z_\lambda$, $\lambda\in P^+$, defined as in 3.1 and as an algebra
by the
$z_i\in F_q^{\varpi_i}$, $1\leq i\leq n$.
Recall the following theorem of [8, Proposition 6.3, Théorème 
7.3].\medskip\noindent
{\bf Theorem.} {\it There exists an $n\times n$ matrix ${\cal A}$
with coefficients in $\A$ and a column matrix $C_q=(C_{q,i})_{1\leq i\leq
n}$
with $C_{q,i}\in\za$
such that\par\noindent
\item{(i)} ${\cal A}C_q=(z_i-C_{n+1}^i)_{1\leq i\leq n}$, \par\noindent
\item{(ii)} $\phi(C_{q,i})$, $1\le i\leq n$, has degree $i+1$ in $Z$ and they 
generate the polynomial algebra $Z$,\par\noindent
\item{(iii)} ${\cal A}=(a_{i,j}(q-q^{-1})^{j-1})$, with
$a_{i,j}=1^jC_{n+1}^{i-1}-2^jC_{n+1}^{i-2}
+\ldots
+(-1)^{i-1}i^jC_{n+1}^0$\qed\medskip\noindent
}\medskip\noindent
{\bf 4.3.} In the following, we set $U:=U(\sl_{n+1})$.
The aim of this section is to describe the $P^+$-filtration
$(F_\lambda)_{\lambda\in P^+}$
of $U(\sl_{n+1})$ in terms of
the natural $\N$-filtration $(U_k)_{k\in\N}$ of $U(\sl_{n+1})$. \par\noindent
For
$\lambda=\sum_i\lambda_i
\varpi_i\in P^+$, set $z(\lambda)=\sum(i+1)\lambda_i$. 
Set $J:=\phi(\ja)$. Recall that, by Proposition 2.1 (ii), 
$J$ is a minimal primitive ideal.
By Lemma 4.1, we can choose $\hq$ in Theorem 2.4 such that
$$\hq\cap
\bigoplus_{h(\lambda)\leq k} F_q^\lambda\subset
\sum_{0\leq m\leq k} (\hq^{\varpi_1})^m,$$ for all $k$. Set $\ha=\hq\cap\ua$,
as in Theorem 2.4. Set $\ha^\lambda$ and $\za^\lambda$ as in 2.5.
By construction, $\ha\cap(\sum_{h(\lambda)\leq k}F_{q,\lambda})
=\bigoplus _{h(\lambda)\leq k} \ha^\lambda$ and
$\za\cap(\sum_{z(\lambda)\leq k}F_{q,\lambda})
=\bigoplus _{z(\lambda)\leq k} \za^\lambda$.
\medskip\noindent
{\bf Proposition.} {\it Set $U=U(\sl_{n+1})$, $H=\phi(\ha)$,
$H_k=H\cap U_k$, $Z_k=Z\cap U_k$, $k\in\N$, we have
\par\noindent
\item{(i)} $Z\otimes H=U$ via multiplication,\par\noindent
\item{(ii)} $\phi(\sum_{h(\lambda)\leq k}\ha^\lambda
)=H_k$,\par\noindent
\item{(iii)} $\phi(\sum_{z(\lambda)\leq
k}\za^\lambda)=Z_k$,\par\noindent
\smallskip\noindent
Proof.}  (i) follows from Proposition 2.5 (iii), Corollary 2.4 and Theorem 2.2.
\par\noindent
Let's prove (ii). For $k$ in $\N$, set $F_{\A,k}=\fa\cap
(\sum_{h(\lambda)\leq k}F_{q,\lambda})$. Let $\overline\phi$ be the
quotient
morphism $\overline\phi$ : $\fa/\ja\rightarrow \,U/J$.
Then, $\overline\phi$ is surjective. We want to show that
$\overline\phi(\overline{F_{\A,k}})$ is the image $\overline U_k$ of
$U_k$ by the canonical surjection.  First of all,
$\overline{F_q^{\varpi_1}}$
is simple and contains $\overline E_1$. Hence,
$\overline\phi(\overline{F_{\A,\varpi_1}})=
\C\oplus\sl_{n+1}$ and $\overline\phi(\bigoplus_{m\leq k}
\overline{(F_{\A,\varpi_1}})^m)=\overline {U_k}$. By Lemma 4.1, this gives
$$\overline {U_k}\subset \overline\phi(\overline F_{\A,k}).\leqno
(4.3.1)$$
For $\mu$ in $P^+$ and $m$ in $\N$, let $c_\mu^m$ be the coefficient of
$q^m$
in the Kostka-Foulkes fonction $K_{\mu,0}(q)$. For $\lambda\in P^+$,
let $c_{\mu}^\lambda=n_{\lambda,\mu}^0$, be the coefficient of $e^\lambda$
in $Q_\mu$. Then, by the Hesselink formula and the Lascoux-Leclerc-Thibon
theorem,
see 0.3 :
$$\dim\overline {U_k}=\sum_{\mu\in P^+}\sum_{m\leq k}c_\mu^m \dim V(\mu)=
\sum_{\mu\in P^+}\sum_{m\leq k}\sum_{h(\lambda)=m}c_{\mu}^\lambda \dim
V_q(\mu).$$
So, by  1.8, Proposition 2.3 and Theorem 2.4,  $\dim\overline {U_k}=
\dim \overline {F_{q,k}}=\dim(\phi(\overline{F_{\A,k}}))$. Hence, we have
equality
in (4.3.1). This gives (ii) by the hypothesis on $\hq$.
\par\noindent
Let's prove (iii).
Fix $s$, $1\leq s\leq n$. Let $\Delta_s$ be the determinant of the matrix
$(a_{i,j})_{1\leq i,j\leq s}$, see
4.2. It is an easy exercice to use the Van der Monde determinant to
obtain  $$\Delta_s=(-1)^{{s(s-1)\over 2}}
\prod_{m=1}^sm\hbox{!}\not=0.$$
From Theorem 4.2, we obtain that $C_{q,s}=D_{q,s}$ modulo $(q-1)$, for all
$1\leq s\leq n$,
with $D_{q,s}\in (q-q^{-1})^{-(s-1)}\sum_{i=1}^s \Q (z_i-C_{n+1}^i)$. 
The elements
$\phi(D_{q,s})=\phi(C_{q,s})$,
$1\leq s\leq n$, generate the algebra $Z(\g)$. So, by induction on the
$P^+$-filtration
of $\za$,  $D_{q,s}$,
$1\leq s\leq n$, generate the algebra $\za$ by using 
Nakayama's lemma on the finite rank $A$-modules 
$\za\cap F_{q,m}$, $m\in\N$. The left hand term in (ii) is 
generated by
$\prod_sD_{q,s}^{n_s}$, with $\sum(s+1)n_s\leq k$. This implies (ii).
 \qed\smallskip\noindent

Set $H^\lambda=\phi(\ha^\lambda)$, resp. $Z^\lambda=\phi(\za^\lambda)$.
This gives the following theorem :\smallskip\noindent
{\bf Theorem.} {\it The $\N$-natural filtration and the $P^+$-filtration
of $U$, see 2.4, can be compared as follow :
$$U_k=\bigoplus_{h(\lambda)+z(\mu)\leq k} H^\lambda\otimes Z^\mu.$$}\qed  \medskip\noindent
{\bf 4.4.} In this section, we illustrate, in the case $\sl_{n+1}$, some
of
the main
results and techniques of this article.\smallskip\noindent
Examples 1-2-3 are devoted to the cone $\co$.\par\noindent
See [25] for the bijection between the canonical base and semi-standard
Young tableaux
and for the calculation of $\E(b)$, where $b$ is a associated to a given
tableau. See [4, Theorem 4.1] for inequations in Example
3.\smallskip\noindent
{\bf Example 1.} For $\g=\sl_2$, the cone $\co$ is $\N$, where $n\in\N$
corresponds in the crystal base to the semi-standard tableau of
shape $2n\varpi$ and of weight 0.
$$ \vbox{\offinterlineskip
\def\th{\hskip 1mm\vrule height 0.4 pt width 2em}\def\tv{$\vert$}
\def\cc#1{\hfill#1\hfill}\cleartabs
\+\th&\th&\th&\th&\th&\th&\th&\th&\cr
\+\tv\cc{1}&\tv\cc{$\ldots$}&\cc{$\ldots$}&\tv\cc{1}
&\tv\cc{2}&\tv\cc{$\ldots$}&\cc{$\ldots$}&\tv\cc{2}&\tv\cr
\+\th&\th&\th&\th&\th&\th&\th&\th&\cr}$$
{\bf Example 2.} For $\g=\sl_3$, and $w_0=s_1s_2s_1$,
the cone $\co$ is generated as a semi-group
by $(1,0,1)$, $(0,1,0)$, $(0,1,1)$, $(2,0,1)$, corresponding to the
following
tableaux of weight 0.
$$\vbox{\offinterlineskip
\def\th{\hskip 1mm\vrule height 0.4 pt width 2em}\def\tv{$\vert$}
\def\cc#1{\hfill#1\hfill}\cleartabs
\+\th&\th&\cr
\+\tv\cc{1}&\tv\cc{2}&\tv\cr
\+\th&\th&\cr
\+\tv\cc{3}&\tv\cr
\+\th\cr}\hskip 5mm
\vbox{\offinterlineskip
\def\th{\hskip 1mm\vrule height 0.4 pt width 2em}\def\tv{$\vert$}
\def\cc#1{\hfill#1\hfill}\cleartabs
\+\th&\th&\cr
\+\tv\cc{1}&\tv\cc{3}&\tv\cr
\+\th&\th\cr
\+\tv\cc{2}&\tv\cr
\+\th\cr}
\vbox{\offinterlineskip
\def\th{\hskip 1mm\vrule height 0.4 pt width 2em}\def\tv{$\vert$}
\def\cc#1{\hfill#1\hfill}\cleartabs
\+\th&\th&\th&\cr
\+\tv\cc{1}&\tv\cc{2}&\tv\cc{3}&\tv\cr
\+\th&\th&\th&\cr
\+ \cc{}&\cc{}&\cc{}&\cr}\hskip 5 mm
\vbox{\offinterlineskip
\def\th{\hskip 1mm\vrule height 0.4 pt width 2em}\def\tv{$\vert$}
\def\cc#1{\hfill#1\hfill}\cleartabs
\+\th&\th&\th&\cr
\+\tv\cc{1}&\tv\cc{1}&\tv\cc{2}&\tv\cr
\+\th&\th&\th&\cr
\+\tv\cc{2}&\tv\cc{3}&\tv\cc{3}&\tv\cr
\+\th&\th&\th&\cr}$$
In the $\C(q)$ algebra of $\co$, the multiplication corresponds
to concatanation.
We have the relation :
$$\vbox{\offinterlineskip
\def\th{\hskip 1mm\vrule height 0.4 pt width 2em}\def\tv{$\vert$}
\def\cc#1{\hfill#1\hfill}\cleartabs
\+\th&\th&\cr
\+\tv\cc{1}&\tv\cc{2}&\tv\cr
\+\th&\th&\cr
\+\tv\cc{3}&\tv\cr
\+\th\cr}*
\vbox{\offinterlineskip
\def\th{\hskip 1mm\vrule height 0.4 pt width 2em}\def\tv{$\vert$}
\def\cc#1{\hfill#1\hfill}\cleartabs
\+\th&\th&\cr
\+\tv\cc{1}&\tv\cc{2}&\tv\cr
\+\th&\th&\cr
\+\tv\cc{3}&\tv\cr
\+\th\cr}*
\vbox{\offinterlineskip
\def\th{\hskip 1mm\vrule height 0.4 pt width 2em}\def\tv{$\vert$}
\def\cc#1{\hfill#1\hfill}\cleartabs
\+\th&\th&\cr
\+\tv\cc{1}&\tv\cc{3}&\tv\cr
\+\th&\th\cr
\+\tv\cc{2}&\tv\cr
\+\th\cr}=
\vbox{\offinterlineskip
\def\th{\hskip 1mm\vrule height 0.4 pt width 2em}\def\tv{$\vert$}
\def\cc#1{\hfill#1\hfill}\cleartabs
\+\th&\th&\th&\cr
\+\tv\cc{1}&\tv\cc{2}&\tv\cc{3}&\tv\cr
\+\th&\th&\th&\cr
\+ \cc{}&\cc{}&\cc{}&\cr}*
\vbox{\offinterlineskip
\def\th{\hskip 1mm\vrule height 0.4 pt width 2em}\def\tv{$\vert$}
\def\cc#1{\hfill#1\hfill}\cleartabs
\+\th&\th&\th&\cr
\+\tv\cc{1}&\tv\cc{1}&\tv\cc{2}&\tv\cr
\+\th&\th&\th&\cr
\+\tv\cc{2}&\tv\cc{3}&\tv\cc{3}&\tv\cr
\+\th&\th&\th&\cr}$$
In the degenerated algebra :
$$\vbox{\offinterlineskip
\def\th{\hskip 1mm\vrule height 0.4 pt width 2em}\def\tv{$\vert$}
\def\cc#1{\hfill#1\hfill}\cleartabs
\+\th&\th&\th&\cr
\+\tv\cc{1}&\tv\cc{2}&\tv\cc{3}&\tv\cr
\+\th&\th&\th&\cr
\+ \cc{}&\cc{}&\cc{}&\cr}*
\vbox{\offinterlineskip
\def\th{\hskip 1mm\vrule height 0.4 pt width 2em}\def\tv{$\vert$}
\def\cc#1{\hfill#1\hfill}\cleartabs
\+\th&\th&\th&\cr
\+\tv\cc{1}&\tv\cc{1}&\tv\cc{2}&\tv\cr
\+\th&\th&\th&\cr
\+\tv\cc{2}&\tv\cc{3}&\tv\cc{3}&\tv\cr
\+\th&\th&\th&\cr}=0$$
This gives, for the case $\sl_3$, an interpretation of a result of
Joseph-Letzter, [22].\smallskip\noindent
{\bf Example 3.} For $\g=\sl_{n+1}$ and $w_0=s_1\ldots s_ns_1\ldots
s_{n-1}\ldots s_1s_2s_1$. Let
$\alpha_{i,j}:=\alpha_i+\ldots\alpha_j$, $1\leq i\leq j\leq n$, be
the set of positive roots. Then $\co$ is the set of points $(n_{i,j},
1\leq i\leq j\leq n)$ such that
$$\sum_{k=j}^nn_{i,k}-\sum_{k=j+1}^nn_{i+1,k}\leq (\sum_{k\leq k'}n_{k,k'}
\alpha_{k,k'},\alpha_i),\hskip 1cm 1\leq i\leq j\leq n.$$
In Example 4, we calculate $P^+$-exponents of type $\mu$, $\mu\in P^+$, in
a particular case.
Then, we deduce the $\N$-exponents of $\mu$.\smallskip\noindent
{\bf Example 4.} For $\g$ of type $A_3$, let's calculate the
$\N$-exponents
of $\mu=\varpi_1+2\varpi_2+\varpi_3$. By [{\it loc. cit.}, 3.4.2 (i)],
$B(\mu)_0$
is parametrized by the semi-standard Young tableaux of shape $\mu$ with
entries
$(1,1,2,2,3,3,4,4)$. These tableaux are given by
$$\vbox{\offinterlineskip
\def\th{\hskip 1mm\vrule height 0.4 pt width 2em}\def\tv{$\vert$}
\def\cc#1{\hfill#1\hfill}\cleartabs
\+\th&\th&\th&\th&\cr
\+\tv\cc{1}&\tv\cc{1}&\tv\cc{2}&\tv\cc{2}&\tv\cr
\+\th&\th&\th&\th&\cr
\+\tv\cc{3}&\tv\cc{3}&\tv\cc{4}&\tv\cr
\+\th&\th&\th&\cr
\+\tv\cc{4}&\tv\cr
\+\th&\cr}\hskip 5mm
\vbox{\offinterlineskip
\def\th{\hskip 1mm\vrule height 0.4 pt width 2em}\def\tv{$\vert$}
\def\cc#1{\hfill#1\hfill}\cleartabs
\+\th&\th&\th&\th&\cr
\+\tv\cc{1}&\tv\cc{1}&\tv\cc{2}&\tv\cc{3}&\tv\cr
\+\th&\th&\th&\th&\cr
\+\tv\cc{2}&\tv\cc{3}&\tv\cc{4}&\tv\cr
\+\th&\th&\th&\cr
\+\tv\cc{4}&\tv\cr
\+\th&\cr}\hskip 5 mm
\vbox{\offinterlineskip
\def\th{\hskip 1mm\vrule height 0.4 pt width 2em}\def\tv{$\vert$}
\def\cc#1{\hfill#1\hfill}\cleartabs
\+\th&\th&\th&\th&\cr
\+\tv\cc{1}&\tv\cc{1}&\tv\cc{2}&\tv\cc{3}&\tv\cr
\+\th&\th&\th&\th&\cr
\+\tv\cc{2}&\tv\cc{4}&\tv\cc{4}&\tv\cr
\+\th&\th&\th&\cr
\+\tv\cc{3}&\tv\cr
\+\th&\cr}$$ $$
\vbox{\offinterlineskip
\def\th{\hskip 1mm\vrule height 0.4 pt width 2em}\def\tv{$\vert$}
\def\cc#1{\hfill#1\hfill}\cleartabs
\+\th&\th&\th&\th&\cr
\+\tv\cc{1}&\tv\cc{1}&\tv\cc{2}&\tv\cc{4}&\tv\cr
\+\th&\th&\th&\th&\cr
\+\tv\cc{2}&\tv\cc{3}&\tv\cc{3}&\tv\cr
\+\th&\th&\th&\cr
\+\tv\cc{4}&\tv\cr
\+\th&\cr}\hskip 5 mm\vbox{\offinterlineskip
\def\th{\hskip 1mm\vrule height 0.4 pt width 2em}\def\tv{$\vert$}
\def\cc#1{\hfill#1\hfill}\cleartabs
\+\th&\th&\th&\th&\cr
\+\tv\cc{1}&\tv\cc{1}&\tv\cc{2}&\tv\cc{4}&\tv\cr
\+\th&\th&\th&\th&\cr
\+\tv\cc{2}&\tv\cc{3}&\tv\cc{4}&\tv\cr
\+\th&\th&\th&\cr
\+\tv\cc{3}&\tv\cr
\+\th&\cr}\hskip 5mm
\vbox{\offinterlineskip
\def\th{\hskip 1mm\vrule height 0.4 pt width 2em}\def\tv{$\vert$}
\def\cc#1{\hfill#1\hfill}\cleartabs
\+\th&\th&\th&\th&\cr
\+\tv\cc{1}&\tv\cc{1}&\tv\cc{3}&\tv\cc{3}&\tv\cr
\+\th&\th&\th&\th&\cr
\+\tv\cc{2}&\tv\cc{2}&\tv\cc{4}&\tv\cr
\+\th&\th&\th&\cr
\+\tv\cc{4}&\tv\cr
\+\th&\cr}$$ $$
\vbox{\offinterlineskip
\def\th{\hskip 1mm\vrule height 0.4 pt width 2em}\def\tv{$\vert$}
\def\cc#1{\hfill#1\hfill}\cleartabs
\+\th&\th&\th&\th&\cr
\+\tv\cc{1}&\tv\cc{1}&\tv\cc{3}&\tv\cc{4}&\tv\cr
\+\th&\th&\th&\th&\cr
\+\tv\cc{2}&\tv\cc{2}&\tv\cc{4}&\tv\cr
\+\th&\th&\th&\cr
\+\tv\cc{3}&\tv\cr
\+\th&\cr}$$
By [{\it loc. cit.}, 3.4.2 (ii)] and Proposition 1.8, the $P^+$-exponents
of $\mu$ are
respectively $2\varpi_1+\varpi_3$, $\varpi_1+\varpi_2$,
$\varpi_1+\varpi_2+\varpi_3$, $\varpi_1+\varpi_2+\varpi_3$, 
$\varpi_1+2\varpi_3$,
 $2\varpi_2$, $\varpi_2+\varpi_3$. After reordering,
the exponents of $\mu$ are $(3,4,5,5,6,6,7)$.
\smallskip\noindent
Example 5 is devoted to the decription of the Joseph-Letzter filtration
for
$\sl_2$.
\smallskip\noindent
{\bf Example 5.} For $\g=\sl_2$, we have $F_{\varpi}=\C\oplus\g\oplus\C
c$,
where $c$ is the quadratic Casimir element. $F_{n\varpi}=F_{\varpi}^n$.
$F_{n\varpi}$ can be described by its $\n^+$-invariant component
$\bigoplus_{k\leq n}(F^{k\varpi})^{\n^+}=\bigoplus_{0\leq i\leq k\leq n}
\C E^{k-i}c^i$. Moreover,
$U_n=\bigoplus_{0\leq h+z\leq n} H^{h\varpi}\otimes Z^{2z\varpi}$.
\medskip\noindent
 \vskip 3cm\noindent

\def\m{{\goth m}}\def\C{{\math C}}
\centerline{APPENDIX}\vskip 1cm\noindent
{\bf Lemma
 A1.} Let
$M$ be a module on a noetherien local domain $R$,
$\m$ be its maximal ideal, 
$K$ be its fraction field and $k$ be its residue  field. 
Suppose dim$_KK\otimes_R M=$dim$_k
k\otimes_R M=r$. Then, $M$ is a free $R$-module with rank $r$.
\smallskip\noindent
Proof. Let $(b_i)$, $1\leq i\leq r$, be a part of $M$ such that
$(1\otimes b_i)$ 
is a base of the $k$-espace $k\otimes_R M$. Then,  $(b)$ is $A$-free.
By cardinality, it is a base of the 
$K$-space $K\otimes_R M$. Let $N$ be the  $A$-module generated by 
$(b)$. Then, $M/N$ is $\m$-torsion. \par\noindent
Let's prove that $M/N$ has no non zero $\m$-torsion elements. 
Let $m$ in $M$ be the lift of a $\m$-torsion element in $M/N$.
For all  $m$ in $M$, there exists $r$ such that $\m^r m\in N$. 
Choose $r$ minimal non zero. By tensoring by $k$, and using the fact that
each element of $R-\m$ is invertible, we obtain 
that 
$(1\otimes
b)$ is not free, which is absurd. Hence, $r=0$ and this gives the Lemma.
\qed \smallskip\noindent
\smallskip\noindent
{\bf Lemma A2.} Let $K$ be a field and let  $R$ be
a localization of  the $K$-algebra $K[T]$ by a 
 multiplicative part $S$ such that $T\not\in S$. Let $M$ be a $R$-module. 
Suppose that \par\noindent
1) $M[T^{-1}]$ is $R[T^{-1}]$-free with rank $r$, 
\par\noindent
2) $M/TM$ is a  $r$-dimensional 
$R/TR$-space. Then, $M$ is a free $R$-module with rank $r$.
\smallskip\noindent
Proof. First of all, we prove the following assertion.
\smallskip\noindent
{\bf Assertion.} There exist  $e_1, \ldots, e_r$ in $M$ such that $N =
\oplus Re_i$ is free with rank $r$, $\{k\otimes e_i\}$ is a base of
$M/TM$ and $N[T^{-1}]=M[T^{-1}]$. \smallskip\noindent
We proceed by induction. Let $e_1^0, \ldots, e_r^0$ in $M$ such that
$(k\otimes e_i^0)$ 
is a base of the $k$-space $k\otimes_R M$. Set $N_0=\oplus Re_i^0\subset
M$.
If $N_0[T^{-1}]=M[T^{-1}]$, there is nothing to prove. If not, then there
exists $x_1$ in $M$ such that $T^kx_1\not\in N_0$ for all $k$ in $\Z$. Set
$N_1=
N_0+Rx_1$. Then, $N_1$ is a finitely generated torsion free $R$-module,
so, it is a free $R$-module. Localizing by $T$, we obtain that
rk$N_1=r$. There exists a base 
$e_1^1, \ldots, e_r^1$ of the $R$ module $N_1$ and,
as $N_1$ contains $N_0$, $(k\otimes
e_i^1)$ 
is a base of the $k$-space $k\otimes_R M$. Using induction and the fact
that $M[T^{-1}]$ is a
noetherian $R[T^{-1}]$-module gives the assertion.\par\noindent
We want to prove that the module $N$ of the assertion is equal to $M$. It
is sufficient to prove that $N_\m=M_\m$ for all maximal ideal $\m$ of $R$. If $T\not\in
\m$, the equality results from $N[T^{-1}]=M[T^{-1}]$. If $T\in\m$, then
$M_\m/
N_\m$ is $T$-torsion by the previous formula. Moreover, as in the proof of
the previous lemma, $M_\m/
N_\m$ has no nonzero $T$-torsion elements. This ends the proof. \qed
 \smallskip\noindent
{\bf Proposition A4.} Define the polynomial ring  $R_1=\C[X_1,..,X_n]$ 
and the local ring $R_2=\C[T]_{(T)}$. Let 
$K_1=\C(X_1,..,X_n)$ et $K_2=\C(T)$ be their fraction fields.
Set $R=R_1\otimes R_2$ endowed with its natural structure of ring.  
Let $M$ torsion-free $R$-module such that \smallskip\noindent
(i) $M/TM$,  is a free $R/TR$-module with rank  $r$,\par\noindent
(ii)
$K_2M=M_T$ is a free $R_1\otimes 
K_2$-module with rank  $r$.\smallskip\noindent
Then, $M$ is a free $R$-module with rank $r$.
\smallskip\noindent
Proof. From Lemma A2, we deduce that $K_1M$ is a free $K_1R$-module 
with rank $r$. We now have 
to prove that $M$ is finitely generated over $R$. This is
the aim of the next lemma :\smallskip\noindent
{\bf Lemma A3.}  Let $M$ torsion-free $R$-module such that $K_1M$, resp.
$K_2M=M_T$, is a free $K_1\otimes R_2$-module with rank  $r$, resp. 
$R_1\otimes 
K_2$-module with rank  $r$. Then, $M$ is finitely generated as a
$R$-module.\smallskip\noindent
Proof. First remark that $R$ is a UFD domain.
Recall that for all multiplicative part $S$ of $R$, we have
$S^{-1}(\bigwedge^rM)=\bigwedge^r(S^{-1}M)$, as $S^{-1}M$-modules. Let 
$F=\bigwedge^rM$. From the hypothesis, we have that $K_1F$ and $K_2F$
are free with rank 1. Let $e_1^s,\ldots, e_r^s\in M$ such that ${e_i^s}$
is
a base of the $K_sR$-module $K_sM$, $s=1$ or $2$. Set
$e^s:=e_1^s\wedge\ldots\wedge e_r^s\in F$. The element $e^s$ is a base of
the $K_sR$ module $K_sF$. 
Let $x$ be in $F$. We have $x=(\alpha_s/a_s)e^s$ and $e^2=(\gamma/d)e^1$
with
$\alpha,\gamma\in R$, $a_s\in R_s$, $d\in R_1$,
GCD$(\alpha_s,a_s)=$GCD$(\gamma,d)=1$. Hence, we obtain
$\alpha_1/a_1=(\alpha_2\gamma)/(bd)$ and then, 
$r:= \alpha_1d/a_1=\alpha_2\gamma/b\in K_1R\cap K_2R=R$.
In conclusion, $F$ is a subset of $(1/d)e^1$.\par\noindent
Fix an element $y$ in $M$, $y=\sum_i(\mu_i/p_i)e_i^1$,
$\mu_i\in R$, $ p_i\in R_1$, GCD$(\mu_i,p_i)=1$. By the previous
assertion,
there exists $s_i$ such that
$$e_1^1\wedge\ldots\wedge e_{i-1}^1\wedge y\wedge e_{i+1}^1\wedge
\ldots\wedge e_r^1=(\mu_i/p_i)e^1=(s_i/d)e^1.$$
We deduce that $\mu_i/p_i=s_i/d$ and so $y\in(1/d) \sum Re_i$. We obtain
that 
$M$ is finitely generated as a $R$-module.\qed
\smallskip\noindent
{\it end of proof of proposition.} In order to show the proposition, 
it suffices, by the Quillen-Suslin theorem,
 cf. [28, Introduction], to prove that $M$ is a
projective, i.e. locally free, $R$-module. 
Let $\m$ be a maximal ideal of $R$. Let's prove that $M_\m$ is free.
If $T$ does not belong to $\m$, then $M_\m$ is a localization of the
module
$M[T^{-1}]=K_2M$ and  
 so, it is a free module. If $T$ belongs to $\m$, it is sufficient to
prove that 
$M_\m$ is a flat
$R_\m$-module.
 This results from the hypothesis and [16, Corollary 6.9].
\qed

\vskip 1 cm \noindent{\petittitre\centerline{ACKNOWLEDGMENTS}}\vskip 1cm
I thank R. Berger and M. Cretin for helpful discussions on Gr\"obner bases
and M.
Reineke for e-mail conversations on Hall algebras. I am grateful to B.
Leclerc
and S. Zelikson for several informations about crystal bases. I also wish
to thank  J. Germoni for comments and encouragements.
Most of all, I am in debt to P. Baumann for his remarks and corrections.
And last, but not least, the article owes a lot to the careful remarks of
the referee.

\vskip 1 cm \noindent{\petittitre\centerline{BIBLIOGRAPHY}}\vskip 1cm
\parindent=0cm
1. P. BAUMANN. {\it Another proof of Joseph and Letzter's separation
of variables theorem for quantum groups}, preprint.\smallskip
2. T. BECKER and V. WEISPFENNING. {\it Gr\"obner bases : a computational
approach
to commutative algebra}, Graduate texts in Mathematics, 141, New-York,
Springer-Verlag, (1993).\smallskip
3. D.J. BENSON. {\it Representations and cohomology I, Basic
representation
theory
of finite groups and associative algebras}, Cambridge studies in Advanced
Mathematics,
30, Cambridge University Press.\smallskip
4. A. BERENSTEIN and A. ZELEVINSKY. {\it Canonical bases for the quantum
group 
of type $A\sb r$ and
piecewise-linear combinatorics.}, Duke Math. J., 82, (1996), no. 3, 
473-502.\smallskip 
5. A. BERENSTEIN and A. ZELEVINSKY. {\it Tensor product multiplicitie, Canonical bases
and Totally positive varieties}, ArXiv:Math.RT/9912012.\smallskip 
6. N. BOURBAKI. {\it Groupes et Alg\`ebres de Lie, Chap. {\rm VI}},
Masson,
 Paris, (1981).\smallskip
7. R.K. BRYLINSKI. {\it Limits of weight spaces, Lusztig's $q$-analogs,
and
fiberings of adjoint orbits}, J. Amer. Math. Soc., 2, (1989),
517-533.\smallskip
8. P. CALDERO. {\it G\'en\'erateurs du centre de $U_q(sl(N+1))$}, Bull.
Sci. Math., 118, (1994),
177-208.\smallskip
9. P. CALDERO. {\it Elements ad-finis de certains groupes quantiques},
 C. R. Acad.
Sci. Paris, t. 316, Serie I, (1993), 327-329.\smallskip
10. P. CALDERO. {\it Etude des $q$-commutations dans l'alg\`ebre
$U_q(\n)$,}
J. Algebra, 178, (1995), 444-457.\smallskip
11. P. CALDERO. {\it On the $q$-commutations in $U_q(\n)$ at roots of
one},
J. Algebra,
Vol 210, (1998), 557-576.\smallskip
12. C. DE CONCINI and V. G. KAC. {\it Representations of quantum groups at
roots of 1}, Colloque Dixmier, Progress in Math., 92, (1990), 471-506.
\smallskip
13. C. DE CONCINI, V. G. KAC, C. PROCESI. {\it Quantum coadjoint action},
J. Amer. Math. Soc.,
5, (1992), 151-189.\smallskip
14. J. DIXMIER. {\it Alg\`ebres enveloppantes}, Cahiers scientifiques
n$^{\circ}$
37, Gaut\-hi\-ers-Villars, Paris, 1974.\smallskip
15. V. G. DRINFELD. {\it On almost cocommutative Hopf algebras}, Leningrad
 Math. J.,
Vol. I, (1990), n$^{\circ}$ 2,  321-342.\smallskip
16. D. EISENBUD.
{\it Commutative algebra. With a view toward algebraic geometry}, 
Graduate Texts in Mathematics, 150, Berlin: Springer-Verlag. 
xvi, (1995).\smallskip
17. W.H. HESSELINK. {\it Characters of the Nullcone}, Math. Ann., 252,
(1980), 179-182.\smallskip
18. A. JOSEPH. {\it Quantum groups and their primitive ideals,}
 Springer-Verlag, 29, (1995).\smallskip
19. A. JOSEPH. {\it On the Mock Peter-Weyl theorem and the Drinfeld double
of
a double}, Reine Angew. Math., 507 ,(1999), 37-56.\smallskip
20. A. JOSEPH and G. LETZTER. {\it Local finiteness of the adjoint action
for quantized envelopping algebras}, J. Algebra, 153, (1992),
 289-318.\smallskip
21. A. JOSEPH and  G. LETZTER.  {\it Separation of
variables for quantized enveloping algebras}, Amer. J. Math.,
116, (1994), 127-177.\smallskip
22. A. JOSEPH and G. LETZTER.  {\it Verma modules annihilators for
quantized enveloping algebras}, Ann. Sci. Ec. Norm. Sup., 28, N$^0$4,
(1995), 493-526.\smallskip
23. A. KANDRI-RODY and V. WEISPFENNING. {\it
Non-commutative Gr\"obner bases in Algebras of solvable type}, J. Symb.
Comp.,
 9, (1990), 1-26.\smallskip
24. M. KASHIWARA. {\it On Crystal Bases}, Canad. Math. Soc.,
Conference Proceed., 16, (1995), 155-195.\smallskip
25. M. KASHIWARA and T. NAKASHIMA. {\it Crystal graphs for representations
of the
$q$-analogue of classical Lie algebras}, J. Alg., 165, N$^0$ 2, (1994),
295-345.
\smallskip
26. S.M. KHOROSHKIN and V.N. TOLSTOI. {\it The uniqueness theorem
for the Universal $\R$-matrix}, Lett. Math. Phys., 24, n$^O$ 3, (1992),
231-244.\smallskip
27. B. KOSTANT. {\it Lie groups representations on polynomial rings},
Amer. J. Math., 85, (1963), 327-404.\smallskip
28. T.Y. LAM.
{\it Serre's conjecture},
Lecture Notes in Mathematics, 635, 
Berlin - Heidelberg - New York: Springer-Verlag. XV, (1978).\smallskip
29. A. LASCOUX, B. LECLERC and J-Y. THIBON. {\it Crystal graphs and
$q$-analogues of Weight Multiplicities for the Root System A$_n$},
Lett. Math. Phys., 35, (1995), 359-374.\smallskip
30. G. LUSZTIG. {\it Canonical bases arising from quantized enveloping
algebras.
 II}, Prog. Theor. Phys. Suppl., 102, (1990), 175-201.
\smallskip
31.  G. LUSZTIG. {\it Introduction to quantum groups}, Progress in
Mathematics,
 110, Birkauser, (1993).
\smallskip
32. M. REINEKE. {\it Multiplicative Properties of Dual Canonical Bases of
Quantum Groups}, J. Alg., 211, (1999), 134-149.\smallskip
33. C.M. RINGEL. {\it Hall algebras and quantum groups,} Invent. Math.,
101, (1990),
583-592. \smallskip
34. M. ROSSO. {\it Analogues de la forme de Killing et du th\'eor\`eme
 de Harish-Chandra
pour les groupes quantiques}, Ann. Sci. Ec. Norm. Sup., 23, (1990),
445-467.\smallskip
{\it Institut Girard Desargues, UPRES-A-5028  \par
Universit\'e Claude Bernard Lyon I, Bat 101\par
69622 Villeurbanne Cedex, France}\smallskip\noindent
e-mail : caldero@desargues.univ-lyon1.fr
\medskip\noindent 
Keywords : Quantum groups, enveloping algebras, canonical base, harmonic
space.

\end